\documentclass[12pt]{article}
\usepackage{graphicx}
\scrollmode
\usepackage{amsfonts}
\usepackage{amsmath}
\usepackage{amssymb}
\begin{document}
\title{\bf  Coarsening, Nucleation, and the Marked Brownian Web}
\author{L.~R.~G.~Fontes \and M.~Isopi \and C.~M.~Newman \and
K.~Ravishankar}
\maketitle

\newtheorem{defin}{Definition}[section]
\newtheorem{Prop}{Proposition}
\newtheorem{teo}{Theorem}[section]
\newtheorem{ml}{Main Lemma}
\newtheorem{con}{Conjecture}
\newtheorem{cond}{Condition}
\newtheorem{prop}[teo]{Proposition}
\newtheorem{lem}{Lemma}[section]
\newtheorem{rmk}[teo]{Remark}
\newtheorem{cor}{Corollary}[section]
\renewcommand{\theequation}{\thesection .\arabic{equation}}
\newcommand{\beq}{\begin{equation}}
\newcommand{\eeq}{\end{equation}}
\newcommand{\beqn}{\begin{eqnarray}}
\newcommand{\beqnn}{\begin{eqnarray*}}
\newcommand{\eeqn}{\end{eqnarray}}
\newcommand{\eeqnn}{\end{eqnarray*}}
\newcommand{\bprop}{\begin{prop}}
\newcommand{\eprop}{\end{prop}}
\newcommand{\bteo}{\begin{teo}}
\newcommand{\bcor}{\begin{cor}}
\newcommand{\ecor}{\end{cor}}
\newcommand{\bcon}{\begin{con}}
\newcommand{\econ}{\end{con}}
\newcommand{\bcond}{\begin{cond}}
\newcommand{\econd}{\end{cond}}
\newcommand{\eteo}{\end{teo}}
\newcommand{\brm}{\begin{rmk}}
\newcommand{\erm}{\end{rmk}}
\newcommand{\blem}{\begin{lem}}
\newcommand{\elem}{\end{lem}}
\newcommand{\ben}{\begin{enumerate}}
\newcommand{\een}{\end{enumerate}}
\newcommand{\bei}{\begin{itemize}}
\newcommand{\eei}{\end{itemize}}
\newcommand{\bdf}{\begin{defin}}
\newcommand{\edf}{\end{defin}}

\newcommand{\nn}{\nonumber}
\renewcommand{\=}{&=&}
\renewcommand{\>}{&>&}
\renewcommand{\le}{\leq}
\newcommand{\+}{&+&}
\newcommand{\fr}{\frac}
\renewcommand{\r}{{\mathbb R}}
\newcommand{\br}{\bar{\mathbb R}}
\newcommand{\Z}{{\mathbb Z}}
\newcommand{\z}{{\mathbb Z}}
\newcommand{\zd}{\z^d}
\newcommand{\zz}{{\mathbb Z}}
\newcommand{\R}{{\mathbb R}}
\newcommand{\tw}{\tilde{\cal W}}
\newcommand{\E}{{\mathbb E}}
\newcommand{\C}{{\mathbb C}}
\renewcommand{\P}{{\mathbb P}}
\newcommand{\N}{{\mathbb N}}
\newcommand{\var}{{\mathbb V}}
\renewcommand{\S}{{\cal S}}
\newcommand{\T}{{\cal T}}
\newcommand{\W}{{\cal W}}
\newcommand{\X}{{\cal X}}
\newcommand{\Y}{{\cal Y}}
\newcommand{\cm}{{\cal M}}
\newcommand{\cp}{{\cal P}}
\newcommand{\h}{{\cal H}}
\newcommand{\f}{{\cal F}}
\newcommand{\cd}{{\cal D}}
\newcommand{\xt}{X_t}
\renewcommand{\ge}{g^{(\epsilon)}}
\newcommand{\xe}{y^{(\epsilon)}}
\newcommand{\ye}{y^{(\epsilon)}}
\newcommand{\bx}{{\bar y}}
\newcommand{\by}{{\bar y}}
\newcommand{\bxe}{{\bar y}^{(\epsilon)}}
\newcommand{\bye}{{\bar y}^{(\epsilon)}}
\newcommand{\bwe}{{\bar w}^{(\epsilon)}}
\newcommand{\bxz}{{\bar y}}
\newcommand{\bwz}{{\bar w}}
\newcommand{\we}{w^{(\epsilon)}}
\newcommand{\Xe}{Y^{(\epsilon)}}
\newcommand{\Ze}{Z^{(\epsilon)}}
\newcommand{\Ye}{Y^{(\epsilon)}}
\newcommand{\ydo}{Y^{(\d)}_{y_0(\d),s_0(\d)}}
\newcommand{\yo}{Y_{y_0,s_0}}
\newcommand{\tye}{{\tilde Y}^{(\epsilon)}}
\newcommand{\hy}{{\hat Y}}
\newcommand{\ve}{V^{(\epsilon)}}
\newcommand{\Ne}{N^{(\epsilon)}}
\newcommand{\ce}{c^{(\epsilon)}}
\newcommand{\cle}{c^{(\l\epsilon)}}
\newcommand{\xet}{Y^{(\epsilon)}_t}
\newcommand{\hxt}{\hat X_t}
\newcommand{\btn}{\bar\tau_n}
\newcommand{\ct}{{\cal T}}
\newcommand{\rn}{{\cal R}_n}
\newcommand{\nt}{{N}_t}
\newcommand{\lnk}{{\cal L}_{n,k}}
\newcommand{\cl}{{\cal L}}
\newcommand{\bw}{\bar{\cal W}}
\newcommand{\tc}{\tilde{\cal C}_b}
\newcommand{\hxtt}{\hat X_{\ct}}
\newcommand{\txnt}{\tilde X_{\nt}}
\newcommand{\xs}{X_s}
\newcommand{\xn}{\tilde X_n}
\newcommand{\tx}{\tilde X}
\newcommand{\hx}{\hat X}
\newcommand{\txi}{\tilde X_i}
\newcommand{\txij}{\tilde X_{i_j}}
\newcommand{\taxi}{\tau_{\txi}}
\newcommand{\txn}{\tilde X_N}
\newcommand{\xk}{X_K}
\newcommand{\ts}{\tilde S}
\newcommand{\tl}{\tilde\l}
\newcommand{\tg}{\tilde g}
\newcommand{\im}{I^-}
\newcommand{\ip}{I^+}
\newcommand{\hal}{H_\a}
\newcommand{\ba}{B_\a}

\renewcommand{\a}{\alpha}
\renewcommand{\b}{\beta}
\newcommand{\g}{\gamma}
\newcommand{\G}{\Gamma}
\renewcommand{\L}{\Lambda}
\renewcommand{\d}{\delta}
\newcommand{\D}{\Delta}
\newcommand{\e}{\epsilon}
\newcommand{\fes}{\phi^{(\epsilon)}_s}
\newcommand{\fet}{\phi^{(\epsilon)}_t}
\newcommand{\fe}{\phi^{(\epsilon)}}
\newcommand{\pset}{\psi^{(\epsilon)}_t}
\newcommand{\pse}{\psi^{(\epsilon)}}
\renewcommand{\l}{\lambda}
\newcommand{\me}{\mu^{(\epsilon)}}
\newcommand{\re}{\rho^{(\epsilon)}}
\newcommand{\tre}{\tilde{\rho}^{(\epsilon)}}
\newcommand{\nue}{\nu^{(\epsilon)}}
\newcommand{\mbe}{{\bar\mu}^{(\epsilon)}}
\newcommand{\rbe}{{\bar\rho}^{(\epsilon)}}
\newcommand{\mb}{{\bar\mu}}
\newcommand{\rb}{{\bar\rho}}
\newcommand{\mbz}{{\bar\mu}}
\newcommand{\s}{\sigma}
\renewcommand{\o}{\Pi}
\newcommand{\om}{\omega}
\newcommand{\tio}{\tilde\o}
\renewcommand{\sl}{\sigma'}
\newcommand{\si}{\s(i)}
\newcommand{\sit}{\s_t(i)}
\newcommand{\ei}{\eta(i)}
\newcommand{\eit}{\eta_t(i)}
\newcommand{\eot}{\eta_t(0)}
\newcommand{\sil}{\s'_i}
\newcommand{\sj}{\s(j)}
\newcommand{\st}{\s_t}
\newcommand{\so}{\s_0}
\newcommand{\xii}{\xi_i}
\newcommand{\xij}{\xi_j}
\newcommand{\xio}{\xi_0}
\newcommand{\ti}{\tau_i}
\newcommand{\te}{\tau^{(\epsilon)}}
\newcommand{\bt}{\bar\tau}
\newcommand{\tti}{\tilde\tau_i}
\newcommand{\tto}{\tilde\tau_0}
\newcommand{\tei}{T_i}
\newcommand{\ttei}{\tilde T_i}
\newcommand{\tes}{T_S}
\newcommand{\tao}{\tau_0}

\renewcommand{\t}{\tilde t}

\newcommand{\da}{\downarrow}
\newcommand{\ua}{\uparrow}
\newcommand{\ar}{\rightarrow}
\newcommand{\lar}{\leftrightarrow}
\newcommand{\va}{\stackrel{v}{\rightarrow}}
\newcommand{\ppa}{\stackrel{pp}{\rightarrow}}
\newcommand{\dw}{\stackrel{w}{\Rightarrow}}
\newcommand{\Va}{\stackrel{v}{\Rightarrow}}
\newcommand{\Ppa}{\stackrel{pp}{\Rightarrow}}
\newcommand{\la}{\langle}
\newcommand{\ra}{\rangle}
\newcommand{\ep}{\vspace{.5cm}}
\newcommand\sqr{\vcenter{
\hrule height.1mm
\hbox{\vrule width.1mm height2.2mm\kern2.18mm\vrule width.1mm}
\hrule height.1mm}}        

\newcommand{\stack}[2]{\genfrac{}{}{0pt}{3}{#1}{#2}}


\begin{abstract}
Coarsening on a one-dimensional lattice
is described by the voter
model or equivalently by coalescing (or annihilating) random walks
representing the evolving boundaries between regions of constant
color and by backward (in time) coalescing random walks
corresponding to color genealogies.
Asympotics for large time and space on the lattice are
described via a continuum space-time voter model whose boundary motion
is
expressed by the {\it Brownian web\/} (BW) of coalescing
forward Brownian motions. In this paper, we study how
small noise in the voter model, corresponding to the nucleation of
randomly
colored regions, can be treated in the continuum limit.
We present a full construction of
the continuum noisy voter model (CNVM)
as a random {\it quasicoloring\/}
of two-dimensional
space time and derive some of its properties.
Our construction is based on a Poisson
marking of the {\it backward\/}
BW within the {\it double\/} (i.e.,
forward and backward) BW.
\end{abstract}


\section{Introduction}

\setcounter{equation}{0}
\label{sec:intro}

In this paper we construct
the one-dimensional continuum
noisy voter model (CNVM) with $q$ colors (opinions),
including the case $q= \infty$.
This model can be thought of as the scaling limit of the noisy voter
model (NVM) on the one-dimensional lattice $\Z$,
extending some of what has been done for the voter model without
noise~\cite{kn:A2,kn:TW,kn:FINR,kn:FINR1}.
In the present paper we will not address
weak convergence issues about the continuum limit procedure
(in the spirit of~\cite{kn:FINR,kn:FINR1}) but will rather
focus primarily on the continuum model itself (in the spirit
of~\cite{kn:TW}).
  We will however prove convergence of finite dimensional distributions.

In the usual (non-noisy) voter model, the color at each site is
updated after a random exponentially distributed waiting time
(with mean one) by taking on the color of a randomly chosen
neighbor ~\cite{L1,L2}.
For updates in the noisy voter model~\cite{kn:GM},
with probability
$1-p$ the site takes on the color of a randomly chosen neighbor and
with
probability $p$ it takes on a color chosen uniformly
at random from all possible
colors or, in the case $q= \infty$, by taking on a completely
novel color.

 When $q=2$, and the two colors are chosen to be
$=+1$ and $-1$, the noisy voter model
is exactly the stochastic
one-dimensional (nearest-neighbor) Ising model with Glauber (heat-bath)
dynamics, where the noise parameter $p$ is related to the temperature
(with the non-noisy case corresponding to zero temperature),
see,
e.g., \cite{kn:GM,kn:FINS1,kn:FINS2} and Subsection \ref{subsec:finiteq}. For more
background on the non-noisy case in the contexts of zero-temperature
Ising/Potts models and diffusion limited reactions,
 see~\cite{kn:DZ,kn:DHP,kn:FV1}.
The scaling limits for voter models, with or without noise, expressed in
 terms of
Brownian webs with or without marks, should be  the same scaling limits
one
gets for certain stochastic PDE models that arise in a variety of physical
settings,
e.g., those that describe nucleation, diffusion and annihilation of
coherent structures
({\em kinks}) in a regime where they can be regarded as pointwise
objects---see,
e.g., \cite{kn:HL1} (also~\cite{kn:HL})
and references therein. This happens in the stochastic Ginzburg-Landau
equation in the limit of small noise and large damping \cite{kn:FV}
or in a classical
$(1+1)$-dimensional $\phi^4$ field theory at finite (low) temperature
\cite{kn:HL1}.

We remark that there is also interest (see, e.g.,~\cite{kn:FINS1})
in stochastic Ising models where the temperature/noise parameter
is not constant in time (and/or space). These type of voter models
are inhomogeneous in time Markov processes. Although we will not
explicitly consider those types of models in this paper, it is clear
 how to extend
the marking constructions given in Section~\ref{sec:mbw} by using
inhomogeneous in time (and/or space) Poisson marking processes.
Of course, some of the simple formulas for two-point functions
given later in this section become more complicated. We now return
to models that are space-time homogeneous.

In the voter model without noise one naturally expects to see large
blocks of the same color and this indeed happens. With
noise, new colors appear within color blocks.
This introduces two new boundaries (or one,  when the new
color appears exactly at an already existing boundary of two
blocks of different colors).
Most of these new colors and associated boundaries
on the (microscopic) lattice will not survive for very long, but
occasionally some will survive for long enough to nucleate a
macroscopic region
of that new color.
In this paper, we introduce a continuum model which describes the long
time behavior of the lattice
voter model with small noise on the appropriate spatial
scale.

In the voter model without noise, the boundaries evolve
according to coalescing (or annihilating, if the colors on either side
of the coalesced boundary are the same) random walks, which in
the scaling
limit converge to coalescing Brownian motions. The collection of
{\it all\/} coalescing walks converges
to the Brownian web (BW)~\cite{kn:A2,kn:TW,kn:FINR,kn:FINR1}.
In the BW,
space-time (double) points where two distinct paths start correspond
in the continuum voter model
to two {\it possible\/}
boundaries which start at a microscopic distance and survive for
a macroscopic time.
In the continuum version of the noisy voter model, if a new color
appears at
such a double point, then it will survive for
a nonzero (macroscopic) time and that double point is a
nucleation point for the new color. All other
newly created colors can be observed at the microscopic level only, so
they
do not appear in the continuum limit.
Even though most space-time points are not double points and most
double points are not nucleation points, we shall see that the
nucleation points form a dense (countable) subset of the
space-time plane.

Roughly speaking, the CNVM is a coloring of the continuum space-time
plane
using the boundaries emerging from the nucleation
points of the newly
created colors. Two things need to be explained---how the nucleation
points are chosen from among all the double points of the Brownian web,
and how a new color taken on at a nucleation point propagates
to a
region of space-time. We will discuss both of these
in the rest of this section,
and then with more detail
in Section~\ref{sec:mbw} and
Subsection~\ref{subsec:3and4} respectively.

We begin by stating a theorem that describes the nature
of the $q = \infty$ CNVM, followed by one describing
the finite $q$ CNVM. Except near the end of this
section of the paper, we restrict attention,
in both the lattice and continuum settings, to the stationary noisy
voter model with time $t \in (-\infty,\infty)$. A natural object we will
focus upon is the mapping from space-time points to the one or
finitely many nucleation points whose color is eventually inherited by
that space-time point. In the lattice case, one can easily define
 things  so that
the mapping is to a single nucleation point; in the continuum limit it
is more natural to
map onto finitely many points. For the lattice case,
let $\Theta_1^p=(X_1^p,T_1^p)$ denote
that
mapping on $\Z\times\R$ for the  $q = \infty$ voter model on $\Z$
with noise parameter $p$. Using the diffusive scalings $\delta^{-1}x$
and $\delta^{-2} t$, we define, for $\delta > 0$,
the rescaled mapping on $(\delta \Z)\times
(\delta^2 \R)$,
\beqn
\label{eq:Theta}
\Theta_\delta^p (x,t)\, = \, (\delta X_1^p(\delta^{-1}x,\delta^{-2} t),
\delta^2 T_1^p(\delta^{-1}x,\delta^{-2} t)) \, .
\eeqn

\bteo
\label{teo:Thetatheorem}
For each $\lambda > 0$, there is a space-time translation invariant
random mapping $\Theta$ from $\R^2$ to finite subsets of $\R^2$
with the following properties:
\begin{itemize}
\item[(1)] For deterministic $(x,t)$, $\Theta(x,t)$ is almost surely
a singleton and distributed as $(x+B(R_\lambda),t-R_\lambda)$ where
$B(\cdot)$ is a standard Brownian motion and $R_\lambda$ is
exponentially
distributed with mean $1/\lambda$ and is independent of $B(\cdot)$.
\item[(2)] For deterministic $(x_1,t_1),\dots,(x_n,t_n)$:
$\Theta$  evaluated at $(x_1,t_1),\dots,(x_n,$ $t_n)$
is the limit in distribution as $\delta \to 0$
of $\Theta_\delta^{\delta^2 \lambda}$
at $(x_1^\delta,t_1^\delta),\dots,(x_n^\delta,t_n^\delta)$
if $(x_i^\delta,t_i^\delta) \to (x_i,t_i)$ for each $i$.
In addition,
\begin{eqnarray}
\label{eq:connectivity}
\P(\Theta(x_1,t_1)=\dots=\Theta(x_n,t_n))\, = \,
\quad \quad \quad \quad \quad \quad \quad & \, & \\
 = \, \, \lim_{\delta \to 0}\P(\Theta_\delta^{\delta^2
\lambda}(x_1^\delta,t_1^\delta)
=\dots=\Theta_\delta^{\delta^2 \lambda}(x_n^\delta,t_n^\delta))\,.
&\, &
\end{eqnarray}

\item[(3)] Almost surely, the set of all nucleation points,
\beqn
\label{eq:nucleate}
{\cal N} \, \equiv  \, \cup_{(x,t) \in \R^2} \Theta(x,t) \, ,
\eeqn
is a dense countable subset of $\R^2$,
and for each nucleation point $(x',t') \in {\cal N}$, its
color region,
\beqn
\label{eq:coloregion}
C_{(x',t')} \, \equiv \, \{(x,t)\in \R^2 :\, \Theta(x,t) \ni (x',t')\},
\eeqn
is a compact, perfect subset of $\R^2$ with an empty interior.
\item[(4)] Almost surely, the cardinalities $|\Theta(x,t)|$ as $(x,t)$
varies over $\R^2$ take on only the values $1,2$ and $3$. Almost
surely, the
set of all
$(x,t)$ with $|\Theta(x,t)| \geq 2$ has Hausdorff dimension
$3/2$ (and hence zero two-dimensional Lebesgue measure).
\item[(5)]
$\psi(x,t)$, defined as $\P(\Theta(0,0) = \Theta(x,t))$, is expressible
as
\beqn
\label{eq:psi1}
\psi(x,t) \, = \, e^{-\lambda |t|}\E(e^{- \sqrt{2\lambda} |x+B(|t|)|})
\,,
\eeqn
where $B$ is standard Brownian motion.
In particular the autocorrelation function is
\beqn
\label{eq:psi2}
\psi(0,t) \, = \, 
1\, - \, \frac{1}{\sqrt{\pi}}
\int_{-\sqrt{\lambda |t|}}^{\sqrt{\lambda |t|}} \exp(-v^2) dv \, ,
\eeqn
and the equal time correlation function is
\beqn
\label{eq:psi3}
\psi(x,0) \, = \, e^{- \sqrt{2\lambda} |x|} \, .
\eeqn
\item[(6)] For any deterministic time $t$, almost surely: for every
$K<\infty$,
the cardinalities $|\Theta(x,t)|$ as $x$
varies over $[-K,K]$ are all one except for finitely many $x$'s where
the
cardinality is two; for each nucleation point $(x',t')$, the
intersection
of its color region, $C_{(x',t')}$ with the horizontal line
$\R \times \{t\}$ is either empty or consists of finitely many closed
intervals
with nonempty interiors.
\item[(7)] For any deterministic $x$, almost surely: the set of $t$'s
such that
$\Theta(x,t)$ is a singleton has full (one-dimensional) Lebesgue
 measure and
$\Theta(x,t)$ is
a doubleton for all other $t$'s; for every nonempty open interval $I$
of $t$'s,
there are infinitely many distinct color regions, $C_{(x',t')}$, that
intersect
the vertical line segment
$\{x\} \times I$ (including both singleton and doubleton points).
\item[(8)] For $(x',t') \in {\cal N}$, let $C^u_{(x',t')}$
denote the unique-color subset of $C_{(x',t')}$:
\beqn
C^u_{(x',t')}
= \{ (x,t) : \Theta(x,t) = (x',t')\}\,.
\eeqn
Almost surely, for
every $(x',t') \in {\cal N}$, $C^u_{(x',t')}$ has strictly
positive (two-dimensional) Lebesgue measure.
\end{itemize}
\eteo

\brm\label{cantor}
Notice the contrast of the (a.s.) natures of the space time color 
regions $C_{(x',t')}$ on one hand, and the fixed deterministic time
color regions $C_{(x',t')}\cap[\R\times\{t\}]$, on the other hand. The
former are Cantor-like sets, while the latter are (possibly empty)
finite disjoint unions of closed intervals with nonempty interiors.
\erm

To construct the CNVM when
$q$ is finite, i.i.d.~uniformly distributed random variables taking
values in $\{1,\dots,q\}$ are assigned to each of the countably many
nucleation points (of the $q=\infty$ model of
Theorem~\ref{teo:Thetatheorem}).
Composing this random color
assignment with the random
mapping $\Theta$ results in a stationary
random {\it quasicoloring\/} $\Phi(x,t)$ mapping $\R^2$ into finite
subsets
of  $\{1,\dots,q\}$ (with cardinality $1$ or $2$ or $3$). Similary one
constructs the corresponding random quasicoloring $\Phi_\delta^p$ on
$(\delta \Z) \times (\delta^2 \R)$ for the lattice noisy
voter model. As a corollary of
Theorem~\ref{teo:Thetatheorem}, we have the following.
\bteo
\label{teo:Phitheorem} For
$\lambda \in (0,\infty)$ and $q\in\{2,3,\ldots \}$,
there exists a space-time translation invariant random
$\{1,2,\ldots,q\}$-quasicoloring
$\Phi$ of the plane such that its finite dimensional distributions
are the limits
as $\delta \to 0$
of those of $\Phi_\delta^{\delta^2 \lambda}$, the diffusively rescaled
stationary one-dimensional voter model with $q$ colors and noise
parameter
$\delta^2 \lambda$.
In particular, for deterministic $(x,t)$,
and $i,j \, \in \{1,\dots,q\}$,
\beqn
\label{eq:onepoint}
\P(\Phi(x,t)=\{i\})\, = \, q^{-1}\,
\eeqn
and
\beqn
\label{eq:twopoint}
\P(\Phi(0,0)=\{i\},\Phi(x,t)=\{j\})\, = \,
q^{-1}\psi(x,t) \, \delta_{i,j} \, + \, q^{-2}(1-\psi(x,t)),
\eeqn
where $\psi(x,t)$ is $\P(\Theta(0,0)=\Theta(x,t))$,
as given in
Theorem~\ref{teo:Thetatheorem}. For any deterministic time $t$,
almost surely, $\Phi(\cdot,t)$
partitions the line into open single-color intervals and
a locally finite set of points where $\Phi(x,t)$ is a
doubleton $\{i,j\}$ separating two intervals of different colors.
For any deterministic
$x$, there is zero probability that
$\Phi(x,\cdot)$ is constant on any nonempty
open interval of $t$-values.
\eteo

\begin{rmk}
We note that \eqref{eq:onepoint} and \eqref{eq:twopoint} are immediate
consequences
of the coloring procedure in which i.i.d. uniformly random colors are
assigned to each nucleation point. The convergence of the finite
dimensional
color distributions follows from the convergence (given in
Theorem~\ref{teo:Thetatheorem}) to
$\P(\Theta(x_1)=\dots=\Theta(x_n))$ and the standard fact that such
probabilities (like the connectivities in lattice percolation models)
determine algebraically the probabilities of all events involving
partitions of $x_1,\dots,x_m$ by distinct color values.
\end{rmk}

 The single time color distribution of the
(stationary) CNVM (the model whose existence is established in
Theorem~\ref{teo:Phitheorem}) is (of course) an
invariant distribution for the CNVM viewed as a (continuum) spin
dynamics. The above theorems give partial results about the nature
of this distribution. In the case where $q=2$, it can be fully
characterized from the
above and the following facts. In this case, the noisy voter model on
the lattice corresponds to
a stochastic Ising model (as we will discuss in more detail in
Subsection~\ref{subsec:finiteq}). Hence, its invariant distribution
is the Gibbs distribution for a nearest neighbour Ising model on $\z$,
which is simply a stationary two-state (spatial) Markov chain.
It is natural to expect the analogous facts to be valid in the
continuum limit and the following theorem states that this is
indeed the case.

\bteo
\label{teo:q=2}
In the special case when $q=2$, for $\lambda \in (0, \infty)$
and deterministic $t$, $\Phi(\cdot, t)$ partitions the line into
single color intervals whose lengths are independent and exponentially
distributed with mean $\sqrt{ 2/\lambda}$. I.e., $\Phi(x,t)$
as a random function of the continuous variable $x$, for deterministic
$t$, is the stationary
two-state
Markov chain with transition rate $\sqrt{\lambda/2}$ from each
state to the other.
\eteo

\begin{rmk} For $2 < q < \infty$, the continuum limit of a nearest
neighbor $q$-state Potts model on $\z$ is a stationary $q$-state
        (spatial)
Markov
chain $\Psi(x)$ for $x \in \r$ with transition rate $r/(q-1)$ from each
state to any of the other $q-1$ states. This process does not appear to
   agree
with the fixed time $t$ CNVM coloring process $\Phi(x,t)$ for any
$q>2$;
that claim can be verified at least for large $q$ as follows. To have
agreement
of the two-point functions, one must take $r = \sqrt{2\lambda} (q-1)/q$,
but then it can be shown that as $\varepsilon \to 0$,
\beqn
\label{eq:gibbs}
\P(\Psi(-\varepsilon)=\Psi(+\varepsilon) \neq \Psi(0))\, \approx \,
\frac{r^2}{q(q-1)^2}\varepsilon^2\, = \,
 \frac{2\lambda}{q^3}\varepsilon^2,
\eeqn
while
\begin{eqnarray}
\label{eq:3point}
\P(\Phi(-\varepsilon,0)=\Phi(+\varepsilon,0) \neq \Phi(0,0))
\quad \quad \quad \quad  \quad \quad & \, & \\
\quad \quad \quad \quad  \geq  \,
\frac{q-1}{q} \P(\Theta(-\varepsilon,0)=\Theta(+\varepsilon,0)
\neq \Theta(0,0)) & \, & \\
\quad \quad     \approx  \, \frac{q-1}{q}C \lambda \varepsilon^2 \, ,
\quad \quad \quad \quad & \, &
\end{eqnarray}
for a universal constant $C>0$. These two formulas do not agree for
sufficiently large $q$.
\end{rmk}

To give even a preliminary construction of the
continuum nucleation mapping $\Theta$ (and hence of the set
${\cal N}$ of nucleation points and their color regions $C_{(x',t')}$),
we need to review some of the properties of the Brownian web and its
associated dual web (\cite{kn:A2,kn:TW,kn:FINR,kn:FINR1,kn:FINR2}).
The BW is a random collection of
paths with specified starting points in space-time. For determinisic
starting points, there are almost surely unique paths starting from
those points and they are distributed as coalescing standard
(except for their starting point) Brownian
motions. But there are random (realization-dependent) double and
triple points from which
respectively two and three paths start.

\begin{figure}[!ht]
\begin{center}
\includegraphics[width=10cm]{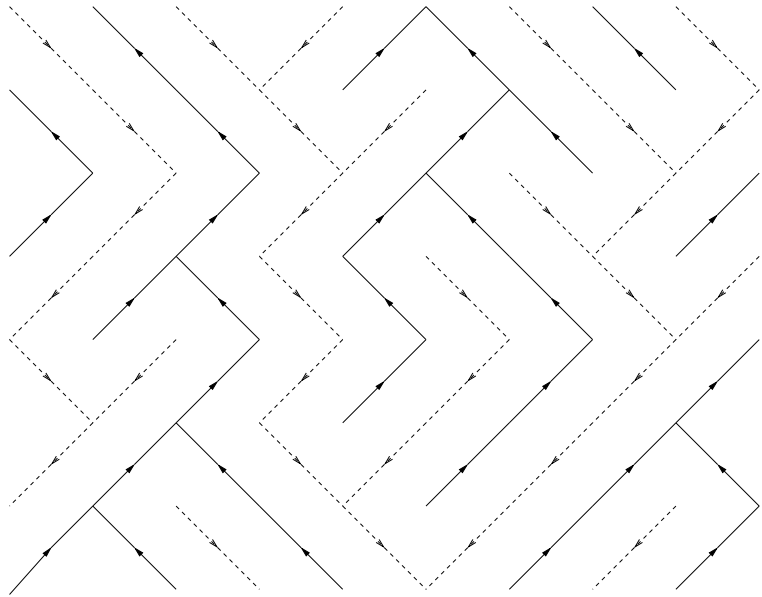}
\caption{Forward coalescing random walks (full lines) in discrete time
and their
dual backward walks (dashed lines).}\label{figdouble}
\label{double} 
\end{center}
\end{figure}

Coexisting with the Brownian web
paths are backward paths; Figure 1 shows the corresponding forward and
backward walks in discrete space and discrete time (vertical
 coordinate).
The forward and backward
Brownian web paths have the following properties
(\cite{kn:TW,kn:FINR1,kn:FINR2}): (1) from
a determinstic point $(x,t)$ with unique forward path $W_{(x,t)}$,
the backward path is the locus of space-time points separating
starting points whose forward path coalesces with $W_{(x,t)}$ from the
left vis-a-vis starting points whose forward path coalesces with
$W_{(x,t)}$ from the
right, (2) every point which is {\it passed through\/} by a backwards
path (i.e., is not merely a starting point for that backwards path) is
either a double point or a triple point (for forward paths),
and (3) the distribution
of all backward paths is exactly that of a time-reversed BW.

Now we can begin to explain the nucleation mapping $\Theta$.
A natural procedure is to have Poisson processes along paths of
the {\em dual} Brownian web of backward paths, which we for now think
of as a BW with
paths going backward in time and ``reflecting'' on the (forward) BW
paths. As already mentioned,
the paths in the dual BW are the loci of the double points of
the forward web.
Follow a (Brownian) path in the dual web and mark it according to
a Poisson process in time with intensity $\lambda$. We do this
for every path in the dual BW
in such a way that on every path segment the markings have intensity
$\lambda$.
This selects a random countable
dense set of double points (of the forward BW).
Almost surely, marks will only appear
on points {\it passed through\/} by a dual path and hence only on
(some of the) double points of the forward web.
We will give a more precise construction in the
next section of the paper.
$\Theta(x,t)$ is now defined by taking
the first marked point backwards in time along every path of the dual
 BW
starting from $(x,t)$. (In order to satisfy the claim of
Property (3) in Theorem~\ref{teo:Thetatheorem}
that the color regions are closed,
we additionally include
$(x,t)$ itself if it is a nucleation point.) In the remaining
sections of the paper,
we will explore this construction in more detail.

We conclude this section by briefly turning to the
case of non-stationary voter models. In this out of equilibrium
setting,
a fairly arbitrary assignment of colors to the points
on the horizontal line $\R \times {t_0}$ (with $t_0$ and the
color assignment deterministic) will lead to a random
quasi-coloring of $\R \times (t_0,\infty)$.
As we have seen, the
nucleation points in $\Theta(x,t)$ are simply the locations of
the first marks (going backwards in time) on the one or more paths
in the backward BW starting from $(x,t)$.
In the non-stationary setting, the corresponding
$\Theta_{t_0}(x,t)$ includes those nucleation marked points $(x',t')$
only if
$t' > t_0$; if the first mark $(x',t')$ along a backward path has
$t'< t_0$ (we leave out $t'=t_0$ since that is a zero probability
event),
then $(x',t')$ is replaced in $\Theta_{t_0}(x,t)$ by the location
$(x'',t_0)$ where the backward path crosses the horizontal line at
$t_0$.
Whether or not $q$ is finite, the colors assigned to these $(x'',t_0)$
points are simply those given by the deterministic initial conditions;
the colors assigned to the nucleation points $(x',t')$ with
$t'>t_0$ are random, as in the stationary setting.
We remark that almost surely, only countably many such $(x'',t_0)$'s
arise as $(x,t)$ varies over $\R \times (t_0,\infty)$.

The rest of the paper is organized as follows. Section~\ref{sec:mbw}
reviews some aspects of the BW and then
presents in detail the Poisson marking process of the BW.
Section~\ref{sec:double} gives background (and some new results)
about the double BW (i.e., the forward BW jointly with the backward BW).
Section~\ref{sec:Proof} presents the proofs
of Theorems \ref{teo:Thetatheorem}, \ref{teo:Phitheorem},
and~\ref{teo:q=2}, including (in Subsection~\ref{subsec:3and4})
providing more details about the color regions
$C_{(x',t')}$ of the CNVM.
There is also an appendix, which is used for one part of
Section~\ref{sec:double} (the analysis of the Hausdorff
dimension of the set of type (1, 2)
points of the BW).


\section{Marked Brownian Web}
\setcounter{equation}{0}
\label{sec:mbw}

In this section we construct
the {\em marked Brownian web} (MBW) as a collection of coalescing
marked Brownian paths. In our application
to the continuum voter model to be discussed later on, we actually
mark the backward (in time) paths.
Before explaining the markings, we review some features of the
(unmarked)
Brownian web. As in~\cite{kn:FINR, kn:FINR1}, we use three metric
spaces:
$(\br^2,\rho)$, $(\o,d)$ and $(\h,d_\h)$.
The elements of the three spaces are respectively: points
in space-time, paths with specified starting points in space-time and
collections of paths with specified starting points.
The BW will be an $(\h,\f_\h)$-valued random variable,
where $\f_\h$ is the Borel $\s$-field  associated to the metric $d_\h$.
Complete definitions of the three metric spaces are given at the end of
this section.
The next theorem, taken from~\cite{kn:FINR1}, gives some of the key
properties of the BW.

\bteo
\label{teo:char}
There is an \( ({\cal H},{\cal F}_{{\cal H}}) \)-valued random variable
\(
\bar{\W} \)
whose distribution is uniquely determined by the following three
properties.
\begin{itemize}
 \item[(o)]  from any deterministic point \( (x,t) \) in $\r^{2}$,
 there is almost surely a unique path \( {W}_{x,t} \)
starting from \( (x,t) \).

 \item[(i)]  for any deterministic \( n, (x_{1}, t_{1}), \ldots,
(x_{n}, t_{n}) \), the joint distribution of \(
{W}_{x_{1},t_{1}}, \ldots, {W}_{x_{n},t_{n}} \) is that
of coalescing Brownian motions (with unit diffusion
 constant), and

\item[(ii)]  for any deterministic, dense countable subset
 \( {\cal   D} \) of \( \r^{2} \), almost surely, \( \bar{\W} \) is the
closure in
\( ({\cal H}, d_{{\cal H}}) \) of \( \{ {W}_{x,t}: (x,t)\in {\cal D} \}. \)
\end{itemize}
\eteo

In our marking procedure the only points in the plane that will be
marked
are those points $(x,t)$ such that a BW path from some time $t' < t$
passes though $(x,t)$.  As previously noted, throughout this
section we will be marking the forward BW,
but later when we deal with the noisy voter model, we will then
work with the marked dual (backward in time) BW.

For each point $(x,t)$, we define the {\it age\/} $\tau(x,t)$ of
that point as the supremum of the set
\[
\{s : \mbox{there exists a path passing through
$(x,t)$ from time $t-s$}\} \, .
\]
All marked points will have {\it strictly\/} positive age.

We proceed with the presentation of four different but
(distributionally) equivalent constructions of the MBW.

\bigskip
\noindent {\bf Construction Via Age-Truncation}

We start by defining the $\varepsilon$-{\it age-truncation\/} of the BW
for any $\varepsilon>0$ as follows.
For each realization of the BW, consider the set ${\cal T}_\varepsilon$
of all points $(x, t)$ in the
plane with age $\tau(x, t)> \varepsilon$.
Next shorten every path in the web by removing (if necessary)
the initial segment consisting of those points of age $\tau \leq \e$.
${\cal T}_\varepsilon$ is the union of the graphs of all these
age-truncated
BW paths and it is almost surely ``locally sparse,'' in the sense
that for
every bounded set $U$, the intersection
${\cal T}_\varepsilon \cap U$ equals the intersection of $U$
with the union of {\it finitely\/} many continuous path segments
(which may be chosen to be disjoint). The locally sparse property
can be verified as follows:
it is known, see \cite{kn:A2}, that for any
any $t$, the intersection
${\cal T}_\varepsilon\cap(\r\times\{t\})$
is (almost surely) locally finite for all $\e>0$. By intersecting
${\cal T}_\varepsilon$ with horizontal strips of height (in
the time variable) $\varepsilon/2$, one sees that there are only
locally finitely many paths passing through the strip.

We now mark each (disjoint) path segment in ${\cal T}_\varepsilon$
according to a
Poisson process in time with rate $\lambda$.
Consider now a sequence of $\varepsilon$'s decreasing to 0.
The marking procedures described above can be carried out for each
$\varepsilon$ and can be coupled in an obvious way so that
the marking for the whole sequence of positive $\varepsilon$'s can be
realized on the same probability space.
Taking the union over all positive $\varepsilon$'s of
these markings gives our first  construction of the marked BW.

Given any $\varepsilon>0$, we denote by ${\cal M}_\varepsilon$ the set
of all marked points in the space-time plane. {\it Conditional\/} on
the BW realization, and hence on the set ${\cal T}_\varepsilon$ (the
``trace'' of the $\varepsilon$-age-truncated BW),
${\cal M}_\varepsilon$
is a spatial Poisson process on the plane with intensity measure
$\lambda \, \mu_\varepsilon$, where $\mu_\varepsilon$ is the locally
finite measure that assigns to each age-truncated path segment in
${\cal T}_\varepsilon$, its $t$-coordinate Lebesgue measure.
The main drawback of this construction is that for any
(bounded) subset $U$ of the plane with nonempty interior,
$\mu_\varepsilon(U) \to \infty$ as $\varepsilon \to 0$
(this is proved in Subsection~\ref{subsec:3and4}) so that
$\lim_{\varepsilon \to 0}\mu_\varepsilon$ is unpleasant to
deal with as a measure on $\r^2$. Our next construction remedies that
feature by using the age as a third coordinate.

\bigskip
\noindent {\bf Construction Via 3D Embedding}

The set ${\cal T}_\varepsilon$ of all $(x,t)$ with age
   $\tau(x,t)>\varepsilon$
is a tree graph embedded continuously in $\r^2$. In our second
construction,
we lift ${\cal T}_\varepsilon$ into $\r^3$ (or more accurately, into
$\r^2 \times (0,\infty)$)
so that we may let $\varepsilon \to 0$ and still
have a locally sparse set. However, the resulting 3D set,
\beqn
{\cal T}^3 = \{(x,t,\tau):\, \tau=\tau(x,t)>0\}\, ,
\eeqn
is no longer a
connected tree graph, but rather
consists of disconnected segments of curves.
The projection of each segment onto
the $(x,t)$-plane is a segment of a path in the BW. Notice that each
segment
in 3D ends when its 2D projection coalesces with another segment that
has an earlier starting point, so that the age of the point of
coalescence is strictly greater than than limit of the age
as the segment (that is about to stop) approaches the coalescence
point. (This age-based priority rule for stopping or
continuing at points of coalescence underlies our next construction.)
We remark that it is natural to regard the 3D curve segments as
being relatively
open, i.e., they do not include either the starting ($\tau=0$) or
ending point.

We may now define a measure $\mu^3$ on $\r^2 \times (0,\infty)$ which
is supported on ${\cal T}^3$ and
assigns to each curve segment its $t$-coordinate
Lebesgue measure. We also define ${\cal M}^3$ as the spatial
Poisson process on $\r^2 \times (0,\infty)$ whose intensity
measure is $\lambda \, \mu^3$. Note that $\mu^3(A)<\infty$ (almost
surely
with respect to the realization of the BW) for
any bounded {\it closed\/} $A$ contained in
$\r^2 \times (0,\infty)$.
The projection of ${\cal M}^3$ onto the $(x,t)$-plane is our
random collection of marked points which
is a.s. countable and dense in $\r^2$.
Note that every point $(x',t')$ in ${\cal M}^3$ has $\tau(x,t)>0$.

In our next construction, we explain how the $2D$ projection of the
connected components of ${\cal T}^3$ may be defined direclty in $\r^2$
without recourse to a $3D$ embedding.

\bigskip
\noindent {\bf Construction Via Tip-Path  Correspondence}

To follow the construction we are about to present, some working
knowledge of the double (forward jointly with dual backward) BW is
needed; this may be obtained by first skimming Section~\ref{sec:double}
and Subsection~\ref{subsec:3and4}.
We are going to mark the paths of the forward BW, taking
advantage of the backward BW, using
properties that hold almost surely.

\begin{figure}[!ht]
\begin{center}
\vspace{-2cm}
\includegraphics[width=17cm]{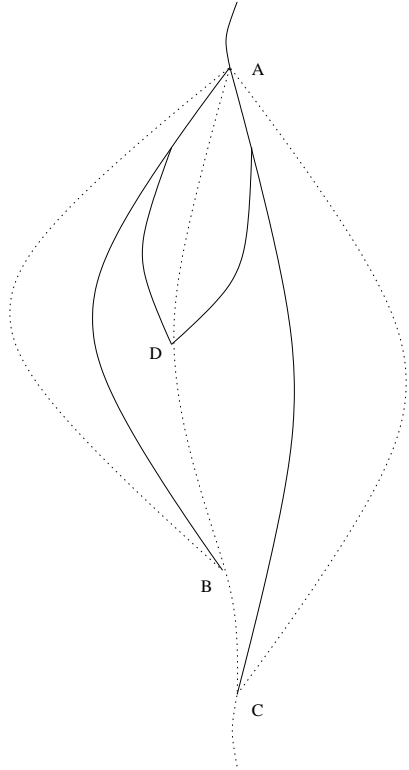}
\vspace{-3cm}
\caption{Forward (upward) paths are in full lines; backward (downward)
ones are dotted. A is a coalescence point of the forward web, and a
starting
tip of the backward bubbles formed by the four
(overlapping) backward path segments,
namely two from A to B, and two from A to C.
The forward path associated to A is the forward
one from B to A. D is a double point of the forward web, and a starting
 tip
of a forward bubble (see Section~\ref{subsec:3and4}), for which A is
an ending tip.}\label{figbubbles}
\label{bubbles} 
\vspace{-.5cm}
\end{center}
\end{figure}

Each coalescence (type (2,1)) point of the {\it forward\/} web is the
starting tip of two {\it backward\/}
bubbles (with disjoint interiors) -- see Subsection~\ref{subsec:3and4}
and Figure~\ref{bubbles}.
For each such point, associate
the subpath of the forward web starting at the highest ending tip of
the two ending tips (one for each of the two backward bubbles), staying within
the respective bubble, and ending at the coalescence point/starting tip of the
respective bubble -- see Figure~\ref{bubbles}. To be consistent
with the 3D embedding construction, this subpath should
be taken relatively open at both ends.

This one-to-one association of subpaths to coalescence points yields a
countable
family (because the set of all coalescence points (of the forward web)
are countable) of disjoint BW path segments. Every
subpath is the initial segment of a path of the BW
belonging to the countable family of all the middle paths starting at
all the triple (type $(0, 3)$) points (each ending tip of a backward
bubble is a triple point, and the chosen
subpath starting there, since it is required to stay within the bubble,
is an initial segment of the middle path from that triple point).
Every triple point will play this role with the initial segment of
its middle path ending when that path coalesces at a point of larger age.
The marking may now be done by using independent rate-$\lambda$
Poisson processes
in time, one for each of the countably many segments.

\bigskip

The three constructions we have presented thus far all use the notion
of the age $\tau(x,t)$ either explicitly, or implicitly in the tip-path
correspondence construction where an age-based precedence relation
between coalescing paths of the BW determines which segment
continues past the coalescence point.
In the tip-path
correspondence construction, one chooses a particular ``skeleton''
(as in (ii) of Theorem~\ref{teo:char}) in which the initial points
of the skeleton are not from a {\it deterministic\/} dense countable
set ${\cal D}$, but rather are the triple
points of the BW realization. We proceed to present a construction
in which one can use any deterministic ${\cal D}$. The main awkward
feature of that construction is that it is not a priori clear
that the resulting MBW has a distribution not depending on the
choice of  ${\cal D}$.

\bigskip
\noindent {\bf Sequential Construction}

This construction begins with independent Brownian
paths starting from any deterministic dense countable
subset ${\cal D}$ of $\R^2$. Mark each Brownian path with marks
from a Poisson process of rate $\lambda$.
One way to do this marking is to consider the set
of Brownian paths ${\cal W}:=\{W_{x,t},\,(x,t)\in{\cal
D}\}$, where $W_{x,t}$ denotes the path starting at
$(x,t)$, and an independent i.i.d.~family ${\cal
N}:=\{N_{x,t},\,(x,t)\in{\cal D}\}$ of Poisson process
(in the time coordinate)
of rate $\lambda$. Now mark the path $W_{x,t}=(f(s),s)_{s\geq
t}$ at the points $(f(S_i),S_i)_{i\geq1}$, where $S_1,S_2,\ldots$
are the successive event times of $N_{x,t}$ after $t$.
Let us denote the marked path thus obtained $W^\ast_{x,t}$ and the
set of marked paths ${\cal W}^\ast:=\{W_{x,t}^\ast,\,(x,t)\in{\cal
D}\}$. Now introduce the set of coalescing marked paths
$\tilde{\cal W}^\ast:=\{\tilde W_{x,t}^\ast,\,(x,t)\in{\cal D}\}$,
as in~\cite{kn:FINR1}, by imposing a precedence relation on the
set of marked paths (note that this is {\it not\/} the
precedence relation based on age used previously but a simpler
one just based on some initial deterministic ordering of ${\cal D}$.
The first coalescing marked path of
$\tilde{\cal W}^\ast$ is the first marked path of ${\cal W}^\ast$.
The $(n+1)$-st coalescing marked path of $\tilde{\cal W}^\ast$ is
formed first with the portion of the $(n+1)$-st marked path of
${\cal W}^\ast$ until it first hits any of the $n$ first
coalescing marked paths of $\tilde{\cal W}^\ast$; from then on, it
follows that marked path (the one it has first hit).

The Brownian web $\bar{\W}$ (as in Theorem~\ref{teo:char}~(ii))
is the closure of the paths in $\tilde{\cal W}^\ast$. It is
important to note however, that {\em there are no new marks
in $\bar{\W}$ beyond those already in the marked skeleton
$\tilde{\cal W}^\ast$}.


\begin{rmk}
When this procedure is used for marking the dual BW then each mark is a
double
point of the forward web which is a starting point of a bubble.
The total time that that bubble exists from its initial to its final
point is identical
to the age of the dual web marked point that coincides with the
forward web double
point (i.e., the age in the dual web equals the bubble lifetime in
the forward web).
\end{rmk}

\bigskip

We end this section with the precise definition of our three metric
spaces.
$(\br^2,\rho)$ is the completion (or compactification) of $\R^2$ under
the
metric $\rho$, where
\begin{equation}
\label{rho}
\rho((x_1,t_1),(x_2,t_2))=
\left|\frac{\tanh(x_1)}{1+|t_1|}-\frac{\tanh(x_2)}{1+|t_2|}\right|
\vee|\tanh(t_1)-\tanh(t_2)|.
\end{equation}
$\br^2$ may be thought as
the image of $[-\infty,\infty]\times[-\infty,\infty]$
under the mapping
\begin{equation}
\label{compactify}
(x,t)\leadsto(\Phi(x,t),\Psi(t))
\equiv\left(\frac{\tanh(x)}{1+|t|},\tanh(t)\right).
\end{equation}

For $t_0\in[-\infty,\infty]$, let $C[t_0]$ denote the set of functions
$f$ from $[t_0,\infty]$ to $[-\infty,\infty]$ such that $\Phi(f(t),t)$
is continuous. Then define
\begin{equation}
\o=\bigcup_{t_0\in[-\infty,\infty]}C[t_0]\times\{t_0\},
\end{equation}
where $(f,t_0)\in\o$ represents a path in $\br^2$ starting at
$(f(t_0),t_0)$.
For$(f,t_0)$ in $\o$, we denote by $\hat f$ the function that extends
$f$
to all
$[-\infty,\infty]$ by setting it equal to $f(t_0)$ for $t<t_0$. Then we
take
\begin{equation}
\label{d}
d((f_1,t_1),(f_2,t_2))=
(\sup_t|\Phi(\hat{f_1}(t),t)-\Phi(\hat{f_2}(t),t)|)
\vee|\Psi(t_1)-\Psi(t_2)|.
\end{equation}
$(\o,d)$ is a complete separable metric space.

Let now $\h$ denote the set of compact
subsets of $(\o,d)$, with $d_\h$ the induced Hausdorff metric, i.e.,
\begin{equation}
\label{dh}
d_\h(K_1,K_2)=\sup_{g_1\in K_1}\inf_{g_2\in K_2}d(g_1,g_2)\vee
\sup_{g_2\in K_2}\inf_{g_1\in K_1}d(g_1,g_2).
\end{equation}
$(\h,d_\h)$ is also a complete separable metric space.

Our description of the continuum noisy voter model, of which the
continuum stochastic Ising model is a particular case, will
involve apart from the MBW, also the dual web to the unmarked BW.
In the next section, we describe the dual BW and the joint object,
the double BW.


\section{Dual and Double Brownian Webs}

\setcounter{equation}{0}
\label{sec:double}

In this section, we construct and characterize the {\em double
Brownian web}, which combines the Brownian web with a {\em dual
Brownian web} of coalescing Brownian motions moving backwards in
time.

In the graphical representation of Harris
for the one-di\-men\-sio\-nal voter model~\cite{kn:H}, coalescing
random walks forward in time and coalescing dual random walks
backward in time (with forward and backward walks not crossing
each other) are constructed simultaneously (see, e.g., the
discussion in~\cite{kn:FINS1,kn:FINS2}). Figure~\ref{figdouble}
provides an example in discrete time. Note that there is no
crossing between forward and backward walks --- a property that is
shared holds also for the {\em double Brownian web} (DBW), which
can be seen as their scaling limit. The
simultaneous construction of forward and (dual) backward Brownian
motions was emphasized in~\cite{kn:TW,kn:STW} and their approach
and results can be applied to extend both the characterization and
convergence results of~\cite{kn:FINR1} to the DBW which includes
simultaneously the forward BW and its dual backward BW.

Our construction and analysis of its properties will rely on a
paper~\cite{kn:STW} of Soucaliuc, T\'oth and Werner together with
results from~\cite{kn:FINR, kn:FINR1} on the (forward) Brownian web
(see also ~\cite{kn:FINR2}).

We begin with an (ordered) dense countable set ${\cal
D}\subset\r^2$, and a family of i.i.d.~standard B.M.'s
$B_1,B_1^b,B_2,B_2^b,\ldots$ and construct
forward and backward paths
$W_1,W_1^b,W_2,W_2^b,\ldots$ starting from $(x_j,t_j)\in{\cal D}$:
\begin{eqnarray}
   \label{eq:for}
   W_j(t)&=&x_j+B_j(t-t_j),\,t\geq t_j\\
   \label{eq:bac}
   W_j^b(t)&=&x_j+B_j^b(t_j-t),\,t\leq t_j.
\end{eqnarray}
Then we construct coalescing and ``reflecting'' paths
$\tilde W_1,\tilde W_1^b,\ldots$ inductively, as follows.
  \begin{eqnarray}
   \label{eq:tw1}
   \tilde W_1&=&W_1;\quad\tilde W_1^b\,\,=\,\,W_1^b;\\
   \label{eq:twn}
   \tilde W_n&=&
   CR(W_n;\tilde W_1,\tilde W_1^b,\ldots,\tilde W_{n-1},\tilde
W_{n-1}^b);\\
   \label{eq:twnd}
   \tilde W_n^b&=&
   CR(W_n^b;\tilde W_1,\tilde W_1^b,\ldots,\tilde W_{n-1},\tilde
W_{n-1}^b),
\end{eqnarray}
where the operation $CR$ is defined in~\cite{kn:STW}, Subsubsection
3.1.4.
We proceed to explain $CR$ for the simplest case, in the definition of
$\tilde W_2$.

As pointed out in~\cite{kn:STW}, the nature of the reflection of a
forward Brownian path $\tilde W$ off a backward Brownian path $\tilde
W^b$
(or vice-versa) is
special. It is actually better described as a push of $\tilde W$ off
$\tilde W^b$
(see Subsection 2.1 in~\cite{kn:STW}). It does not have an explicit
formula in general, but in the case of one forward path and one
backward
path, the form is as follows. Following our notation and construction,
we ignore $\tilde W_1$ and consider $\tilde W_1^b$ and  $\tilde W_2$
in the time interval $[t_2,t_1]$ (we suppose $t_2<t_1$; otherwise,
$\tilde W_1^b$ and  $\tilde W_2$ are independent). Given $W_2$ and
$\tilde W_1^b$,
for $t_2\leq t\leq t_1$,
\begin{equation}
   \label{eq:refl}
   \tilde W_2(t)=
   \begin{cases}
     W_2(t)+\sup_{t_2\leq s\leq t}(W_2(s)-\tilde W_1^b(s))^-,
           \mbox{ if } W_2(t_2)>\tilde W_1^b(t_2);\\
     W_2(t)-\sup_{t_2\leq s\leq t}(W_2(s)-\tilde W_1^b(s))^+,
           \mbox{ if } W_2(t_2)<\tilde W_1^b(t_2).
   \end{cases}
\end{equation}
After $t_1$, $\tilde W_2$ interacts only with $\tilde W_1$, by
coalescence.

We call
${\cal W}^D_n:=\{\tilde W_1,\tilde W_1^b,\ldots,\tilde W_n,\tilde
W_n^b\}$
{\em coalescing/reflecting forward and backward Brownian motions
(starting at
$\{(x_1,t_1),\ldots,(x_n,t_n)\}$)}.
We will also use the alternative notation ${\cal W}^D({\cal D}_n)$ in
place of
${\cal W}^D_n$, where ${\cal D}_n:=\{(x_1,t_1),\ldots,(x_n,t_n)\}$.
\brm
\label{rm:stw}
In Theorem 8 of~\cite{kn:STW}, it is proved that the above
construction is a.s.~well-defined, gives a perfectly
coalescing/reflecting system (see Subsubsection 3.1.1
in~\cite{kn:STW}),
and for every $n\geq1$, the distribution of ${\cal W}^D_n$
does not depend on the ordering of ${\cal D}_n$.
It also follows from that result that $\{\tilde W_1,\ldots,\tilde
W_n\}$ and
$\{\tilde W_1^b,\ldots,\tilde W_n^b\}$ are separately forward and
backward coalescing
Brownian motions, respectively. Thus
$\{\tilde W_1,\tilde W_2,\ldots\}$ and
$\{\tilde W_1^b,\tilde W_2^b,\ldots\}$ are forward and backward
Brownian web skeletons, respectively.
\erm

\brm
One can alternatively use a set ${\cal D}^b$ of starting points for the
backward
paths different than ${\cal D}$
rather than our choice above of ${\cal D}^b={\cal D}$.
\erm

We now define dual spaces of paths going backward in time $(\Pi^b,d^b)$
and
a corresponding $({\cal H}^b,d_{{\cal H}^b})$
in an obvious way, so that they are the dual versions of
$(\Pi,d)$ and $({\cal H},d_{{\cal H}})$,
and then define ${\cal H}^D={\cal H}\times{\cal H}^b$ and
$$d_{{\cal H}^D}((K_1,K_1^b),(K_2,K_2^b))=\max(d_{{\cal
H}}(K_1,K_2),d_{{\cal H}^b}(K_1^b,K_2^b)).$$

As in the construction of the (forward) BW, we now define
\beqn
{\cal W}_n^D({\cal D})
&=&
\{\tilde W_1,\ldots,\tilde W_n\}\times\{\tilde W_1^b,\ldots,\tilde
W_n^b\},\\
{\cal W}^D({\cal D})
&=&\{\tilde W_1,\tilde W_2,\ldots\}\times\{\tilde W_1^b,\tilde
W_2^b,\ldots\},\\
\bar{\cal W}^D({\cal D})&=&
\overline{\{\tilde W_1,\tilde W_2,\ldots\}}\times\overline{\{\tilde
W_1^b,\tilde W_2^b,\ldots\}}.
\eeqn
The latter closures are in $\Pi$ for the first factor and in $\Pi^b$
for the second one.

{}From Remark~\ref{rm:stw}, we have that
$$\bar{\cal W}:=\overline{\{\tilde W_1,\tilde W_2,\ldots\}} \mbox{ and }
\bar{\cal W}^b:=\overline{\{\tilde W_1^b,\tilde W_2^b,\ldots\}}$$
are forward and backward Brownian webs, respectively.
We proceed to state three propositions and one theorem; their proofs
follow directly from the results
and methods of~\cite{kn:STW, kn:FINR, kn:FINR1}

\bprop
\label{prop:dbwcpt} Almost surely, $\bar{\cal W}^D({\cal
D})\in{\cal H}^D$  (i.e. $\overline{\{\tilde W_1,\tilde
W_2,\ldots\}}$ and $\overline{\{\tilde W_1^b,\tilde
W_2^b,\ldots\}}$ are compact).
\eprop

\brm
\label{rm:dlim} It is immediate from this proposition that
$$\bar{\cal W}^D({\cal D})=\lim_{n\to\infty}{\cal W}_n^D({\cal D}),$$
\erm
where the limit is in the $d_{{\cal H}^D}$ metric.

\bprop
\label{prop:dpaths}
$\bar{\cal W}^D({\cal D})$ satisfies
\begin{itemize}
\item[$(o^D)$] From any deterministic $(x,t)$ there is almost surely
a unique forward path and unique backward path.
\item[$(i^D)$] For any deterministic
${\cal D}'_n:=\{(y_1,s_1),\ldots,(y_n,s_n)\}$ the
forward and backward paths
from ${\cal D}'_n$, denoted $\bar{\cal W}^D({\cal D},{\cal D}'_n)$,
are distributed as coalescing/reflecting forward and backward Brownian
motions starting at ${\cal D}'_n$. In other words,
$\bar{\cal W}^D({\cal D},{\cal D}'_n)$ has the same distribution as
${\cal W}^D({\cal D}'_n)$.
\end{itemize}
\eprop

\bprop
\label{prop:dind}
The distribution of $\bar{\cal W}^D({\cal D})$ as an
$({\cal H}^D,{\cal F}_{{\cal H}^D})$-valued
random variable (where
${\cal F}_{{\cal H}^D}={\cal F}_{{\cal H}}\times{\cal F}_{{\cal
H}^b}$),
does not depend on ${\cal D}$.
Furthermore,
\begin{itemize}
\item[$({ii}^{D})$] for any deterministic dense ${\cal D}'$, almost
surely
\begin{equation*}
\bar{\cal W}^D({\cal D})=
\overline{\{W_{x,t}:(x,t)\in{\cal
D}'\}}\times\overline{\{W_{x,t}^b:(x,t)\in{\cal D}'\}},
\end{equation*}
\end{itemize}
where $W_{x,t},W_{x,t}^b$ are respectively the forward and backward
paths in $\bar{\cal W}^D({\cal D})$
starting from $(x,t)$, and
the closures in $({ii}^{D})$ are in $\Pi$ for the first factor and in
$\Pi^b$ for the second one.
\eprop

\bteo
\label{teo:chard}
The double Brownian web is characterized
(in distribution, on $({\cal H}^D, {\cal F}_{{\cal H}^D})$) by
conditions $(o^D),(i^D)$ and $(ii^D)$.
\eteo

We now discuss ``types'' of points $(x,t)\in\r^2$,
whether deterministic or not.
For the (forward) Brownian web,
we define
\beqn\nn
m_{\mbox{{\scriptsize in}}}(x_0,t_0)\!\!&=&\!\!
\lim_{\varepsilon\downarrow0}\{\mbox{number of paths in } \W
\mbox{ starting at some }
t_0-\varepsilon
\mbox{ that pass}\\
\label{eq:min}
&& \mbox{through } (x_0,t_0) \mbox{ and are disjoint for }
t_0-\varepsilon<t<t_0\};\\\nn
m_{\mbox{{\scriptsize out}}}(x_0,t_0)\!\!&=&\!\!
\lim_{\varepsilon\downarrow0}\{\mbox{number of paths in } \W
\mbox{ starting at }
(x_0,t_0)\mbox{ that are}\\
\label{eq:mout}
&& \mbox{ disjoint for }
t_0<t<t_0+\varepsilon\}.
\eeqn

For $\W^b$, we similarly define $m^b_{\mbox{{\scriptsize
in}}}(x_0,t_0)$ and
$m^b_{\mbox{{\scriptsize out}}}(x_0,t_0)$.

\begin{figure}[!ht]
\begin{center}
\includegraphics[width=6cm]{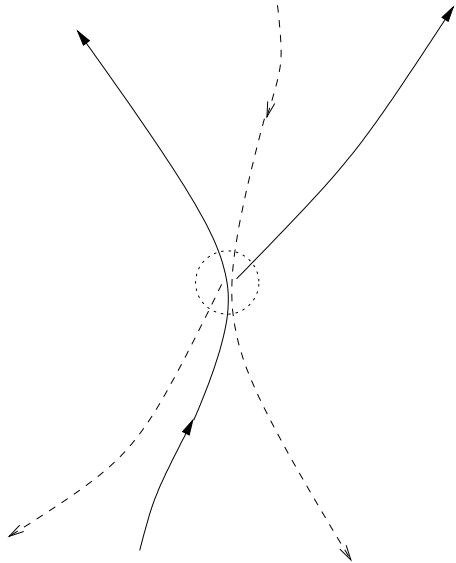}
\caption{A schematic diagram of a point $(x_0,t_0)$ of type
$(m_{\mbox{{\scriptsize\emph{in}}}},m_{\mbox{{\scriptsize
\emph{out}}}})=(1,2)$,
with necessarily also
$(m^b_{\mbox{{\scriptsize\emph{in}}}},m^b_{\mbox{{\scriptsize
\emph{out}}}})=(1,2)$.
In this example the incoming forward path connects to the leftmost
outgoing path
(with a corresponding dual connectivity for the backward paths);
at some of the other points of type $(1,2)$ it
will connect to the rightmost path.}\label{fig12}
\label{12} 
\end{center}
\end{figure}

\bdf
The type of $(x_0,t_0)$ is the pair $(m_{\mbox{{\scriptsize
\emph{in}}}},
m_{\mbox{{\scriptsize \emph{out}}}})$---see Figure~\ref{fig12}.
We denote by $S_{i,j}$ the set of points of $\r^2$ that are of type
$(i,j)$, and by $\bar S_{i,j}$ the set of points of $\r^2$ that are of
type $(k,l)$ with $k\geq i$, $l\geq j$.
\edf

\brm
\label{rm:dense}
Using the translation and scale invariance properties of the
Brownian web distribution, it can be shown that
for any $i,j$, whenever $S_{i,j}$ is nonempty, it must be dense in
$\r^2$.
The same can be
said of $S_{i,j}\cap\r\times\{t\}$ for deterministic $t$.
These denseness properties can also be shown for each $i,j$ by more
direct arguments.
\erm

\bprop
\label{prop:dual}
For the double Brownian web, almost surely for {\em every}
$(x_0,$ $t_0)$ in $\r^2$, $m^b_{\mbox{{\scriptsize \emph{in}}}}
(x_0,t_0)=m_{\mbox{{\scriptsize \emph{out}}}}(x_0,t_0)-1$
and
$m^b_{\mbox{{\scriptsize \emph{out}}}}(x_0,t_0)=
m_{\mbox{{\scriptsize \emph{in}}}}(x_0,t_0)+1$. See Figure~\ref{fig12}.
\eprop

\noindent{\bf Proof } It is enough to prove
(i) that for every incoming forward path to a point $(x,t)$,
there are two locally disjoint backward paths starting at that point
with one on either side of the forward path;
and (ii) that for every two locally disjoint backward paths
starting at a point $(x,t)$, there is an
incoming forward path to $(x,t)$ between the two backward paths.
(Note that by a $t \longleftrightarrow -t$ time reflection
argument, one would then get a similar result for incoming backward
paths and pairs of outgoing forward paths.)

Let us start with the first assertion. Let $\g$ be an incoming
forward path to $(x,t)$. This means that the starting time $s$ of
$\g$ is such that $s<t$. By Proposition 4.3 of~\cite{kn:FINR1},
the portion of $\g$ above time $s+\varepsilon$ is in the
forward skeleton for every $\varepsilon>0$. Now consider a sequence of
pairs of backward paths $(\g_k,\g'_k)$ starting at
$((x_k,t_k),(y_k,s_k))\in{\cal D}\times{\cal D}$ with
$((x_k,t_k),(y_k,s_k))\to((x,t),(x,t))$ as $k\to\infty$,
$s+\varepsilon <
s_k,t_k<t$, $x_k<\g(t_k)$ and $y_k>\g(s_k)$. From the reflection
of the forward and backward skeletons off each other and the fact
that two backward paths in the skeleton must coalesce once they
meet, it follows that $\g_k(t')<\g'_k(t')$ for all
$t'\in[\mbox{max}\{s_k,t_k\},t']$. We then conclude from
compactness that there are two locally disjoint limit paths, one
for $(\g_k)$ and one for $(\g'_k)$, both starting from $(x,t)$.

We argue (ii) similarly. Given two locally disjoint backward paths
$\g,\g'$ starting at $(x,t)$, there exists $s<t$ such that
either $\g(t')<\g'(t')$ for $s<t'<t$ or
$\g'(t')<\g(t')$ for $s<t'<t$. Suppose it is the first case;
otherwise, switch labels. Then choose a point $(x',s')\in{\cal D}$
with $s<s'<t$ and $\g(s')<x'<\g'(s')$. The fact that the portions of
$\g$ and $\g'$ below time $t-\varepsilon$
is in the backward skeleton for every $\varepsilon>0$ and the
reflection
of the forward and backward skeletons off each other now implies that
the forward path starting at $(x',s')$ is squeezed between
$\g$ and $\g'$ and goes to $(x,t)$.

\bteo
\label{teo:types}
For the (double) Brownian web, almost surely, every $(x,t)$ has one of
the
following types, all of which occur: $(0,1)$, $(0,2)$, $(0,3)$,
$(1,1)$,
$(1,2)$, $(2,1)$.
\eteo

\brm
Points of type $(1,2)$ are particularly
interesting in that the single incident path continues along exactly
one
of the two outward paths --- with the choice determined intrinsically
rather than by some convention. See Figure~\ref{fig12} for a schematic
diagram of a ``left-handed'' continuation. An $(x_0,t_0)$ is of
type $(1,2)$ precisely if both a forward and a backward path pass
through
$(x_0,t_0)$. It is either left-handed or right-handed according to
whether the forward path is to the left or the right of the backward
path near $(x_0,t_0)$. Both varieties occur and the proof of
Theorem~\ref{teo:typesa} below shows that the Hausdorff dimension
of $1$ applies separately to each of the two varieties.
\erm

T\'oth and Werner~\cite{kn:TW} gave a definition of types of points of
$\r^2$ similar to ours, but for a somewhat different process
and proved the above theorem with that definition and for that process
(see definition at page 385, paragraph of equation~(2.28) and
Proposition 2.4 in~\cite{kn:TW}).
One way then to establish Theorem~\ref{teo:types} is to show the
equivalence
of ours and T\'oth and Werner's definition and that their arguments
hold for our process. We prefer, for the sake of simplicity and
completeness, to give a direct argument, out of which the following
complementary results also follow.

\bteo
\label{teo:typesa}
Almost surely,
$S_{0,1}$ has full Lebesgue measure in $\r^2$,
$S_{1,1}$ and $S_{0,2}$ have Hausdorff dimension $3/2$ each,
$S_{1,2}$ has Hausdorff dimension $1$, and
$S_{2,1}$ and $S_{0,3}$ are both countable and dense in $\r^2$.
\eteo

\bteo
\label{teo:typesb}
Almost surely: for every $t$
\begin{itemize}
\item[a)] $S_{0,1}\cap\r\times\{t\}$ has full Lebesgue measure in
$\r\times\{t\}$;
\item[b)] $S_{1,1}\cap\r\times\{t\}$ and $S_{0,2}\cap\r\times\{t\}$
are both countable and dense in $\r\times\{t\}$;
\item[c)] $S_{1,2}\cap\r\times\{t\}$, $S_{2,1}\cap\r\times\{t\}$ and
$S_{0,3}\cap\r\times\{t\}$ have all cardinality at most $1$.
\end{itemize}

For every deterministic $t$,
$S_{1,2}\cap\r\times\{t\}$, $S_{2,1}\cap\r\times\{t\}$ and
$S_{0,3}\cap\r\times\{t\}$ are almost surely empty.
\eteo

\noindent{\bf Proof of Theorems~\ref{teo:types}
and~\ref{teo:typesa} } We start by ruling out the cases that do
not occur almost surely. For $i,j\geq0$, $S_{i,j}=\emptyset$
almost surely if $j=0$ or $i+j\geq4$. The first case is trivial.
We only need to consider $\bar S_{i,j}$ for the cases $i=3, j=1$
and $i=2, j=2$, since the other ones are either contained or dual
to these. By Proposition 4.3 of~\cite{kn:FINR1}, $\bar
S_{3,1}$ consists of points which are almost surely in the
skeleton and where three paths coalesce. But the event that three
coalescing Brownian paths starting at distinct points coalesce at
the same time is almost surely empty. By Proposition 4.3
of~\cite{kn:FINR}, $\bar S_{2,2}$ consists of points
(almost surely in the double skeleton) where two different forward
paths coalesce and a backward path passes. Since for any two
forward and one backward Brownian paths in the double skeleton,
the event that this happens is almost surely empty, by the
perfectly coalescing/reflecting property of the paths in the
double skeleton (see Subsubsection 3.1.1 and Theorem 8
of~\cite{kn:STW}) the conclusion follows.

Now, for the types that do occur.

{\bf Type (2,1) } By the above, $S_{2,1}=\bar S_{2,1}$ almost
surely, and  $\bar S_{2,1}$ consists almost surely of {\em points
of coalescence}, that is all points where two paths coalesce. By
Proposition 4.3 of~\cite{kn:FINR1}, it is almost surely
a subset of the skeleton, and thus is countable (since there is at
most one coalescence point for each pair of paths starting from
$\cal D$ in the skeleton). It is easy to see that it is dense
since the paths from a pair of nearby points in $\cal D$ also
coalesce nearby with probability close to one.

{\bf Type (1,2) }
By the above, $S_{1,2}=\bar S_{1,2}$ almost surely, and
$\bar S_{1,2}$ consists almost surely of points where forward paths
meet
backward paths. Thus, it is a subset of the (union of the traces of all
the paths in the)
skeleton. It is easy to see that it is almost surely nonempty (and also
dense).
We need only consider two such paths, say $W$ and $W^b$,
the former a forward one starting at $(0,0)$ (without loss of
generality, by the
translation invariance of the law of $\W^D$), and the latter a backward
one starting at an
arbitrary deterministic $(x_0,t_0)$, with $t_0>0$ to avoid a trivial
case. It is clear that
the random set $\Lambda$ of space-time points $(t,W(t))$ for times
$t\in[0,t_0]$ when
$W(t)=W^b(t)$ has a positive, less than
one probability of being empty. We will argue next the following claim.

{\noindent\bf Claim}
{\it $\Lambda$ has Hausdorff dimension $1$ for almost every
  pair of trajectories $(W,W^b)$
for which it is nonempty.}

By Proposition~\ref{prop:dpaths}, the distribution of
$\{(W(t),W^b(t)):\,0\leq t\leq t_0\}$
(which is all that matters for this) can be
described in terms of two (forward) independent standard Brownian
motions $B,B^b$ as follows
(see equations~(\ref{eq:for})-(\ref{eq:refl})). Let
$W^b(t)=x_0+B^b(t_0-t),\,t\leq t_0$, and
$\tau=\inf\{t\in[0,t_0]:\,B(t)=W^b(t)\}$,
with $\inf\emptyset=\infty$. If $\tau=\infty$, then $W=B$; otherwise,
$W(t)=B(t)$
for $0\leq t\leq\tau$, and for $\tau\leq t\leq t_0$,
\beqnn
  W(t)=
   \begin{cases}
B(t)+\sup_{0\leq s\leq t}(W^b(s)-B(s)),
           \mbox{ if } W^b(0)<0;\\
B(t)-\sup_{0\leq s\leq t}(B(s)-W^b(s)),
           \mbox{ if } W^b(0)>0.
   \end{cases}
\eeqnn
Rewriting in terms of $W'(t):=W^b(t)-W^b(0),\,0\leq t\leq t_0$, which
is a standard
Brownian motion independent of $B$, we have (for $0\leq t\leq t_0$)
$$
   W(t)=
   \begin{cases}
B(t)+\sup_{0\leq s\leq t}\{W'(s)-B(s)\}-W'(t_0)+x_0,
           \mbox{ if } W'(t_0)>x_0;\\
B(t)+\inf_{0\leq s\leq t}\{W'(s)-B(s)\}-W'(t_0)+x_0,
           \mbox{ if } W'(t_0)<x_0,
   \end{cases}
$$
if $\tau\leq t\leq t_0$, with
$\tau=\inf\{t\in[0,t_0]:\,B(t)=W'(t)-W'(t_0)+x_0\}$;
otherwise, $W(t)=B(t)$.

From the above discussion, we conclude that $\Lambda$
has the same distribution as the random set ${\cal G}$ obtained as
follows. Let
${\cal T}^+$ and ${\cal T}^-$ be the sets of positive and negative
record times of the
standard Brownian motion $X(t):=(W'(t)-B(t))/\sqrt2$, respectively,
i.e., ${\cal T}^+$ is the set of $t\geq0$ such that $X(t)=\sup_{0\leq
s\leq t}X(s)$
and ${\cal T}^-$ is the same except with inf in place of sup. Consider
also the
standard Brownian motion $Y(t):=(W'(t)+B(t))/\sqrt2$, which is
independent of $X$.
If $W'(t_0)>x_0$, then
${\cal
G}=\{([(X(t)+Y(t))/\sqrt2]-[(X(t_0)+Y(t_0))\sqrt2]+x_0,t):\,t\in{\cal
T}^+\cap[\tau,t_0]\}$;
if $W'(t_0)<x_0$, then
${\cal
G}=\{([(X(t)+Y(t))/\sqrt2]-[(X(t_0)+Y(t_0))\sqrt2]+x_0,t):\,t\in{\cal
T}^-\cap[\tau,t_0]\}$.

It follows from Proposition~\ref{prop:haus} in Appendix~\ref{app:haus}
that the sets
${\cal G}^\pm:=\{(X(t)+Y(t),t):\,t\in{\cal T}^\pm\cap[0,t_0]\}$ (one
for each sign, respectively)
both have Hausdorff dimension $1$ almost surely.
Since the events $\{W'(t_0)>x_0\}$, $\{W'(t_0)<x_0\}$ and
$\{\tau<t_0\}$ all have
positive probability, the claim follows.

{\bf Type (1,1) } $\bar S_{1,1}$ almost surely consists of {\em
points of continuation} of paths, that is, all points $(x,t)$ such
that there is a path starting earlier than $t$ that touches
$(x,t)$. By Proposition 4.3 of~\cite{kn:FINR1}, $\bar
S_{1,1}$ is almost surely a subset of the skeleton. Since the
trace of any single path has Hausdorff dimension
$3/2$~\cite{kn:Ta} and the countable union of such sets has the
same dimension, it follows that $\bar S_{1,1}$ has Hausdorff
dimension $3/2$ almost surely. By the previous parts of the proof,
$\bar S_{1,1} \setminus S_{1,1}$ has lower dimension and so
$S_{1,1}$ has the same Hausdorff dimension of $3/2$.

{\bf Type (0,1) } We claim that any deterministic point is a.s.~of
this type, hence (by applying Fubini's Theorem) $S_{0,1}$ is
a.s.~of full Lebesgue measure in the plane. That
$m_{\mbox{{\scriptsize in}}}(x_0,t_0)=0$ a.s.~for every
deterministic $(x_0,t_0)$ follows from Proposition 4.3
of~\cite{kn:FINR1}, since, if $m_{\mbox{{\scriptsize
in}}}(x_0,t_0)\geq1$, then there would be a path in the skeleton
passing through $(x_0,t_0)$, but this event clearly has
probability zero. The assertion that $m_{\mbox{{\scriptsize
out}}}(x_0,t_0)=1$ a.s.~for every deterministic $(x_0,t_0)$ is
property $(o)$ of Theorem \ref{teo:char}.

By Proposition~\ref{prop:dual}, the remaining types $(0,2)$ and $(0,3)$
are dual respectively to
$(1,1)$ and $(2,1)$, since the other types are dual to these.
Since ${\bar \W}^b$ is distributed like the standard Brownian web
(modulo a time reflection), the claimed results for types $(0,2)$ and
$(0,3)$
follow from those already proved for $(1,1)$ and $(2,1)$.

\noindent{\bf Proof of Theorem~\ref{teo:typesb}}

{\bf Type (0,1) }
$\bar S_{1,1}$ is almost surely in the skeleton, thus making
$\bar S_{1,1}\cap\r\times\{t\}$ countable for all $t$. By a duality
argument,
the same is true for $\bar S_{0,2}$. Since $\bar S_{0,1}=\r^2$ a.s.~by
Theorem~\ref{teo:types}, it follows that a.s.~for
all $t$, $S_{0,1}\cap\r\times\{t\}$ is of full Lebesgue
measure in the line.

Again, of the remaining types, it is enough by duality to
consider $(1,1)$, $(2,1)$ and $(1,2)$.

{\bf Type (2,1) }
For any deterministic $t$ and $(x_i,t_i)$ with $t_i<t$, $i=1,2$,
the probability that two independent Brownian paths starting
at $(x_i,t_i)$, $i=1,2$, respectively, coalesce exactly at time $t$
is zero. Since $S_{2,1}$ is in the skeleton,
$S_{2,1}\cap\r\times\{t\}=\emptyset$ almost surely.
Now, for any $t$, $|S_{2,1}\cap\r\times\{t\}|>1$ implies that
there are four independent Brownian paths starting at different
points, and such that the coalescence time of the first two and
that of the last two are the same. That this has zero probability
implies that a.s.~for all $t$, $|S_{2,1}\cap\r\times\{t\}|\leq 1$.

{\bf Type (1,2) }
For any deterministic $t$,
$S_{1,2}\cap\r\times\{t\}=\emptyset$ almost surely, since the
probability of two fixed paths, one forward, one backward, meeting
at a given deterministic time is $0$.
Indeed, from the analysis of type $(1,2)$ done above in the proof of
Theorem~\ref{teo:typesa}, this is because the probability that a
Brownian motion has a record value at a given deterministic time
is~$0$.
For any $t$, $|S_{1,2}\cap\r\times\{t\}|>1$ implies that
there exist in the double Brownian web skeleton two pairs,
each consisting of
one forward and one backward path, such that in both pairs
the forward and backward paths meet at the same time.
We claim that this has zero probability
and thus that $|S_{1,2}\cap\r\times\{t\}|\leq 1$ almost surely.
To verify the claim, we again use the analysis of type $(1,2)$ done
for Theorem~\ref{teo:typesa}, which shows that it suffices to
prove that there is zero probability that
two independent standard Brownian motions $B_1,\,B_2$
have a common strictly positive record time. But, as noted
in Appendix~\ref{app:haus}, this is the same as having zero
probability for $B_1,\,B_2$ to both have a zero at
a common strictly positive time. This latter
probability is indeed zero because of
the well known fact that the two-dimensional Brownian motion
$(B_1,B_2)$ a.s.~does not return to $(0,0)$.

{\bf Type (1,1) } Since points with
$m_{\mbox{{\scriptsize in}}}\geq1$ are a.s.~in the skeleton,
$\bar S_{1,1}\cap\r\times\{t\}$ is a.s.~countable
(and easily seen to be dense) for every
$t\in\r$. Now the previous parts of the proof
imply that the same holds for
$S_{1,1}\cap\r\times\{t\}$ for every $t\in\r$.


\section{Proofs of Theorems ~\ref{teo:Thetatheorem},
~\ref{teo:Phitheorem}, and~\ref{teo:q=2}}
\label{sec:Proof}
 \setcounter{equation}{0}
\subsection{Proofs of Theorem ~\ref{teo:Thetatheorem}
(1) and (2) -- convergence}
\label{subsec:1and2}
(1) $\Theta(x,t)$ is defined (see Sections~\ref{sec:intro}
and~\ref{sec:mbw}) by considering {\it all}
paths
in the {\it backward} Brownian web from $(x,t)$ and taking the set of
first marked points (i.e. closest in time to $t$) of those paths. The
marking is done with rate $\lambda$ and in the special case where
$(x,t)$
is itself a mark, $\Theta(x,t)$ includes both $(x,t)$ and the first
mark
$(x',t')$ with $t' < t$. Property (1)
~\ref{teo:Thetatheorem} follows
from Property~(o) of  Theorem ~\ref{teo:char}.

(2) We recall that $\Theta_\delta^{\delta^2 \lambda}(x,t)$ is the
value of $\Theta$ for a rescaled process where time is scaled like
$\delta^2$, space is scaled like $\delta$ and the nucleation rate
is $\delta^2 \lambda$. Therefore it follows from well known
results that the coalescing random walks starting from
$(x_i^\delta,t_i^\delta), 1 \leq i \leq n$ converge in
distribution to coalescing Brownian motions starting from
$(x_i,t_i), 1 \leq i \leq n$. Since the rate of the Poisson clocks
(nucleation rate) is $\delta^2 \lambda$, the markings for the $n$
coalescing random walks for all $\delta$ and the $n$ coalescing
Brownian motions can be done using $n$ fixed Poisson processes
each of rate $\lambda$. That is, the marking part of the process
for all $n$ rescaled walks and the limiting coalescing Brownian
motions can be coupled.

We define the coalescing walks and
Brownian motions by introducing priorities. When a walk (or a BM)
with label $i$ meets a walk (or BM) with label $j <i$ it follows
the path of the walk (or BM) with label j after that time. Let
${T'}_i^\delta, 2 \leq i \leq n$ be the time when the walk starting
from $(x_i^\delta,t_i^\delta)$ meets a walk starting from
$(x_j^\delta,t_j^\delta)$ with $j < i$ and let $T_i^\delta =
 t_i^\delta -
{T'}_i^\delta$ (recall that we are moving backwards in time). Denote by
${T'}_1^\delta = \min ({T'}_2^\delta, \cdots , {T'}_n^\delta)$  the
time when all the $n$ walkers have coalesced and let $T_1^\delta =
t_1^\delta - {T'}_1^\delta$.  Let $T_i, 1 \leq i \leq n$, be the
corresponding times for the Brownian motions starting from
$(x_i,t_i), 1 \leq i \leq n$.

Since the $T_i^\delta$'s are functionals
of the $n$ random walks starting from $(x_i^\delta,t_i^\delta)$,
it follows that not only the walks, but also
$T_i^\delta, 1 \leq i \leq n$, converge in joint
distribution to the continuum paths and $T_i, 1 \leq i \leq n$.
Property (2) then follows; for example, to
prove the second claim of Property (2), we
observe that $\P(\Theta_\delta^{\delta^2
\lambda}(x_1^\delta,t_1^\delta) =\dots=\Theta_\delta^{\delta^2
\lambda}(x_n^\delta,t_n^\delta))=\E e^{- \lambda \sum_1^n
T_i^\delta}$. Now since $T_i^\delta, 1 \leq i \leq n$ converge in
distribution to $T_i, 1 \leq i \leq n$ we have

\begin{eqnarray*}
 \P(\Theta(x_1,t_1)=\dots=\Theta(x_n,t_n)) &=& \E (e^{- \lambda
\sum_1^n
T_i}) \\
& = & \lim_{\delta \to 0} \E (e^{- \lambda \sum_1^n T_i^\delta})\\
&= & \lim_{\delta \to 0}\P(\Theta_\delta^{\delta^2
\lambda}(x_1^\delta,t_1^\delta) =\dots=\Theta_\delta^{\delta^2
\lambda}(x_n^\delta,t_n^\delta))\, .
\end{eqnarray*}

\subsection{Proofs of Theorem ~\ref{teo:Thetatheorem} (3) and
(4) -- nucleation points and color regions}
\label{subsec:3and4}

The MBW, in the context of Theorem ~\ref{teo:Thetatheorem}, is a
dual/backwards dynamics, in the
sense that it is the continuum version of
marked coalescing random walks, which is dual to the noisy voter model
and runs backwards in time. In this
way we get an indirect/dual/backwards description of the continuum
version of the noisy voter model (CNVM).  We can
get a direct/forward description of the CNVM by
considering not only the
forward web, but also simultaneously the dual web. The dual web is
needed in order to get the marks placed on the
dual paths.

 Once the marks are in place, we can
focus on the paths of the forward web
starting at the marks in $\r^2$. We note that, since the marks are on
(non-starter points of) dual paths,
each one is a double point of the forward web,
and thus is the origin of a ``bubble" (of
the color it was assigned). We note that bubbles will occur inside
other bubbles,
with the color of the inside one prevailing.

\begin{figure}[!ht]
\begin{center}
\includegraphics[width=3cm]{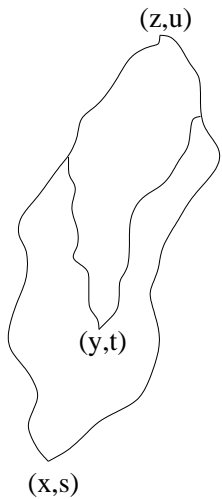}
\caption{\small Two space-time bubbles with starting tips at $(x,s)$
and $(y,t)$, and common ending tip at $(z,u)$.}
\label{fig:clus}
\end{center}
\vspace{-.5cm}
\end{figure}

Let us look at these bubbles, each consisting of the
closed region
of $\r^2$ bounded by the two paths
starting at a mark until they coalesce. In this situation, we
call the mark the {\em starting
tip} of the bubble. These are the nucleation points. We will call the
space time point
where a bubble ends, i.e., the space time point above
the starting tip where the bubble boundary paths meet (and coalesce), the
{\em ending tip} of that bubble. See Figure~\ref{fig:clus}.

We argue that almost surely there are countably many nucleation points,
since almost surely there are countably many marks. This
easily follows from any of the constructions of Section~\ref{sec:mbw};
for example using the sequential construction, one may make the
following two observations:

 i) The marked points occur only on the skeleton of the (dual)
BW and the skeleton is a countable collection of paths (i.e., those
starting at the countable set
${\cal D})$.

ii) On each path of the skeleton the nucleation events (marks) are
Poisson events corresponding to a Poisson process with rate
$\lambda$.

In order to show that almost surely the set of all nucleation
points is dense in $\R^2$ it is sufficient to show that almost
surely for all $\varepsilon>0$
the square of side $\varepsilon$
centered at the origin, $S_0^\varepsilon$, contains
a nucleation point. This can be shown in a variety of ways;
we proceed with one of them. For
large enough $n \in \N$, divide
$S_0^\varepsilon$ into small rectangles of horizontal side length
$n^{-1/4}$ and vertical side length $n^{-1}$. There are
about
$C_1\varepsilon^2 n^{5/4}$ such rectangles in $S_0^\varepsilon$. For
each of those rectangles the probability that there is no path
within the rectangle from the midpoint of the top edge to
the bottom edge is bounded by the probability
of the event that
  a (backward) Brownian motion starting at the midpoint of the top
edge leaves the rectangle through one of the side edges. By
standard arguments, this probability is bounded by $c_1 e^{-c_2
\sqrt n}$. For each of the $n^{-1/4}$ by $n^{-1}$
rectangles, the conditional
probability, given that there is such a path in the rectangle, that
there is no mark on that path in the rectangle is
equal to $e^{- \lambda n^{-1}}$. Since there are $C_1 \varepsilon^2
n^{5/4}$ such rectangles we have that the probability
of no nucleation point in $S_0^\varepsilon$ is bounded above by
\beqn
C \varepsilon^2 n^{5/4} e^{-c_2 n^{1/2}} + e^{- \lambda n^{-1} C_1
\varepsilon^2
n^{5/4}} \to 0 \hbox{ as } n \to \infty.
\eeqn
This proves that almost
surely the set of all nucleation points is dense.

Now we show that $C_{(x',t')}$ for a nucleation point $(x',t')$ is a
compact subset of $\R^2$. That it is bounded is clear since it is
contained in the bubble from $(x',t')$, which is a.s.~bounded. 
It is thus enough to show it is closed. We
start with the space-time bubble starting at $(x',t')$.
This is a closed subset of $\R^2$. Let $(x_1,t_1),
(x_2,t_2), \ldots $ be some ordering of the nucleation
points in the interior of the bubble $B_{(x',t')}$. Then we
claim that
\[
C_{(x',t')} = B_{(x',t')} - \lim_{ n \to \infty}
\hbox{rint}(\cup_{j=1}^n B_{(x_j,t_j)})\\
=  \cap_{n=1}^\infty (B_{(x',t')} -
  \hbox{rint}(\cup_{j=1}^n B_{(x_j,t_j)})
\]
where  rint denotes the relative interior (relative
to $B_{(x',t')}$). This would show that $C_{(x',t')}$ is
closed since it would be the intersection of closed sets.

To justify the claim we note first that $(x_0,t_0) \in
B_{(x',t')}$ is also in $C_{(x',t')}$ if and only if there is a
backwards path from
$(x_0,t_0)$ within $B_{(x',t')}$ which touches none of $(x_1,t_1),
(x_2,t_2), \ldots $. We thus need only show that for all $n$, a point
  $(x_0,t_0) \in B_{(x',t')}$ does {\it not} have a backward
  path from $(x_0,t_0)$ touching any of $(x_1,t_1),
(x_2,t_2), \ldots, (x_n,t_n)$ if and only if $(x_0,t_0)$
belongs to $\hbox{rint} (\cup_{j=1}^n B_{(x_j,t_j)})$.

For a point
$(x_0,t_0)$ in $\hbox{rint} (\cup_{j=1}^n B_{(x_j,t_j)})$, it is
not hard to see that {\it every} backward path from
$(x_0,t_0))$ must pass through one of the $(x_k,t_k),1 \leq k \leq
n$,
(before it reaches $(x',t')$). On the other hand, for
points in $B_{(x',t')}- \hbox{rint} (\cup_{j=1}^n B_{(x_j,t_j)}) $
every backward path enters (in arbitrarily small time)
the interior of $B_{(x',t')}- \hbox{rint} (\cup_{j=1}^n
B_{(x_j,t_j)})
$ and then can no longer touch any of the
$(x_k,t_k), 1 \leq k \leq n$. If it did, the mark at that
$(x_k,t_k)$ would be at a point where two backward paths
coalesced. There is zero probability of such a mark occurring, as can
be seen from, e.g., the sequential construction of the MBW
in Section ~\ref{sec:mbw}.

Perfectness of $C_{(x',t')}$ follows immediately from the
fact that its points belong to nondegenerate (continuous) path
segments which 
are themselves in $C_{(x',t')}$.  That $C_{(x',t')}$ has empty
interior follows immediately from the denseness of ${\cal N}$.

Property (4) is an immediate consequence of
Theorems~\ref{teo:types} and~\ref{teo:typesa}.

\subsection{Proof of Theorem~\ref{teo:Thetatheorem} (5) -- two-point
   functions}
\label{ssec:twop}

The two-point correlation function
$\psi(x,t)$ is defined as $\P(\Theta(0,0) = \Theta(x,t))$.
It is the probability that two (almost surely unique)
backwards paths starting at $(0,0)$ and
$(x,t)$ will not get marks before meeting.

Without loss of generality, we may assume $t \geq 0$.
The case when $t<0$ negative is then easily reduced to this by
time translation-invariance.
Denote by $Y$ the  position at time $0$ of the backward path starting
at $(x,t)$.
Now observe that $\psi(x,t) = \E [\psi(Y,0)e^{-\lambda t}]$
and that
\[
\E[\psi(Y,0)]=\E(e^{-2\lambda \hat T}) \, ,
\]
where $\hat T$ is the time when two independent BMs
starting at $(0, 0$) and $(0,
Y)$ first meet.

By elementary properties of BM, this equals
$\E(e^{-\lambda T})$, where $T$ is the time standard BM first reaches
$|Y|$.
Now this is simply the Laplace transform of the distribution
of a hitting time
of BM.  By the optional sampling theorem (see, e.g., \cite{kn:B}),
it can be proved without calculation (see, e.g., \cite{kn:KT}) that
\[
\E(e^{-\lambda T}) = e^{-\sqrt{2\lambda} |Y|}\, .
\]
So we have
\[
\psi(x,t) = e^{-\lambda t}\E(e^{-\sqrt{2\lambda} |Y|}) \, ,
\]
where $Y$ is distributed as $\mathcal{N}(x,t)$.

Whether or not $ t \geq 0$, we thus have
\[
\psi(x,t) \, = \, e^{-\lambda |t|}\E(e^{- \sqrt{2\lambda} |x+B(|t|)|})
\,,
\]
where $B$ is a standard Brownian motion.
When $t = 0$,
\[
\psi(x,0) = e^{-\sqrt{2\lambda} |x|}
\]
since $Y=x$.
When $x=0$ so that $Y$ is centered, we change variables in the
integral and  get:
\[
\E(e^{-\sqrt{2\lambda} |Y|}) = e^{\lambda
t}\frac{2}{\sqrt{\pi}}\int_{\sqrt{\lambda |t|}}^{\infty} \exp(-y^2) dy
\, ,
\]
yielding that
\[
\psi(0,t)= \frac{2}{\sqrt{\pi}}\int_{\sqrt{\lambda |t|}}^{\infty}
\exp(-y^2) dy.
\]

\subsection{Proof of Theorem~\ref{teo:Thetatheorem} (6) -- fixed
time coloring}
\label{ssec:fincol}

The above descriptions raise a natural question as to how the
color configurations of the MBW dynamics look at fixed
times.
By the direct description we know that we have
space-time color clusters one inside the other almost surely. It
is not difficult to see that each cluster has another cluster
inside it, e.g., by concluding from the scaling of the marked
random walks that the marks of the MBW are dense in $\r^2$. This
might suggest that the latter picture occurs also for fixed
positive times, i.e., the color clusters at positive
times\footnote{By a color cluster at fixed time, say $t$, we mean
any connected component of the intersection with $\r\times\{t\}$
of a space-time color cluster, say $C_0$, minus the intersections
with $\r\times\{t\}$ of the closures of the space-time color
clusters contained in $C_0$.} would also be such that each one has
another cluster inside it. But this is not the case, as one
sees for the case $q=2$ from Theorem~\ref{teo:q=2}.
In fact, we argue next
that, even in the case that $q=\infty$, the color configurations
of the MBW dynamics at fixed positive times have finitely many
clusters in each finite interval almost surely.

It is enough to consider a single-time segment $[a,b]\times\{t\}$
for deterministic $a,b,t$ with $a<b$ and show that the expected
total time-length of all the disjoint pieces of all the backward
paths of the backward BW starting at $[a,b]\times\{t\}$ down to
time $t-s$ is almost surely finite for arbitrary $s>0$.

For $r>0$, let $\eta(t,r;a,b)$ be the number of disjoint points at
$\R\times\{t-r\}$ which are touched by the backward paths starting
on $[a,b]\times\{t\}$. Then the above mentioned total sum can be
expressed as
\begin{equation}
\label{tl1} \int_0^s\eta(t,r;a,b)\,dr.
\end{equation}

We then need to show
\begin{equation}
\label{tl2}
\E\int_0^s\eta(t,r;a,b)\,dr=\int_0^s\E[\eta(t,r;a,b)]\,dr
\end{equation}
is finite. And this follows from the formula
\begin{equation}
\label{tl3} \E[\eta(t,r;a,b)]=(b-a)/\sqrt{\pi r},
\end{equation}
which holds for all $r>0$ (see
Theorem 1.1 in~\cite{kn:FINR}).

\subsection{Proof of Theorem~\ref{teo:Thetatheorem}
(7) -- nonpersistence}
\label{subsec:pers}

We will show in this section that persistence, in its usual sense
of no (or only finitely many) color changes at fixed spatial locations
for
strictly positive amounts of (rescaled) time, does not occur
in the continuum noisy voter model (in contrast to the non-noisy
voter model~\cite{kn:FINS1,kn:DHP}).
More precisely, we show for $q = \infty$ that
\begin{description}
\item {almost surely, any deterministic
vertical interval with nonzero length has infinitely many colors;
the set of points with a unique color has full Lebesgue measure in
the interval; all other points have exactly two colors and there
are infinitely many of them.}
\end{description}
For definiteness, we take $\{0\}\times[0,1]$ as
the deterministic vertical interval. The above claims will follow
from the fact that $L$, the total time-length of all the disjoint
pieces of all the paths of the backward BW starting at
$\{0\}\times[0,1]$ down to time $0$, is almost surely infinite.

We first write $L$ as
\begin{equation}
   \label{eq:inp0}
L=\int_0^1N_s\,ds,
\end{equation}
where, for $0\leq s\leq1$, $N_s$ denotes the number of distinct points
in $\R\times\{1-s\}$ touched by paths starting on $\{0\}\times[1-s,1]$.

We now show that, for $0<s\leq1$, $N_s=\infty$ almost surely.
This implies the above claim.
By rescaling, if $0<s\leq1$, then $N_s$ has the same distribution as
$\hat N_t$, the number of distinct points in $\R\times\{0\}$ touched by
paths starting on $\{0\}\times[0,t]$, for any $t>0$. Now $\hat N_t$ is
nondecreasing in $t$. Let $\hat N_\infty:=\lim_{t\to\infty}\hat N_t$,
the number of distinct points in $\R\times\{0\}$ touched by paths
starting on $\{0\}\times[0,\infty)$.
Then  $N_s$ has the same distribution as $\hat N_\infty$ and it is thus
enough to argue that $\hat N_\infty=\infty$ almost surely.

One straightforward way of arguing the latter point (there are
other slicker arguments that use the Double Brownian Web) is
to show that the event that there
exists a sequence of paths starting on $\{0\}\times[0,\infty)$ which
are
disjoint down to time $0$ has probability 1. For that, it is enough to
exhibit for every $\delta>0$, a sequence $t_1<t_2<\ldots$ such that
the event $A$ that the backward paths starting at
$(0,t_1),(0,t_2),\ldots$
are disjoint down to time $0$ has probability at least $1-\delta$.
The paths from $(0,t_1),(0,t_2),\ldots$ can be taken as independent
Brownian paths.

To choose $t_1,t_2,\ldots$, we start with a sequence $p_1,p_2,\ldots$
such that $p_i>0$ for all $i\geq1$ and $\prod_{i=1}^\infty p_i
   \geq1-\delta$.
We take $t_0=0$ and $M_0=0$ and proceed inductively as follows. Having
defined $t_0,M_0\ldots,t_{n-1},M_{n-1}$, let $t_n$ be such that the
probability of the event $\tilde A_n$ that the path from $(0,t_n)$ does
not touch the rectangle $[-M_{n-1},M_{n-1}]\times[0,t_{n-1}]$ by time
$0$
is at least $(1+p_n)/2$. That there exists such $t_n$ follows from the
fact that for any $t,M>0$, the probability that the path from $(0,t')$
does not touch the rectangle $[-M,M]\times[0,t]$ by time $0$ goes to 1
as $t'\to\infty$. With such a $t_n$ picked, choose $M_n$ such that the
probability of the event $\hat A_n$ that the path from $(0,t_n)$ does
not touch the vertical sides of the rectangle $[-M_n,M_n]\times[0,t_n]$
by time $0$ is at least $(1+p_n)/2$.
That there exists such $M_n$ follows from the
fact that for any $t>0$, the probability that the path from $(0,t)$
does not touch the vertical sides of the rectangle $[-M,M]\times[0,t]$
by time $0$ goes to 1 as $M\to\infty$.

Now let $A_n=\tilde A_n\cap\hat A_n$. Then
$A\supset\cap_{n=1}^\infty A_n$,
and
\begin{equation}
   \label{eq:prob}
\P(A)\, \geq \, \prod_{n=1}^\infty \P(A_n) \, \geq \,
\prod_{n=1}^\infty [(1+p_n)/2]^2 \, \geq \, \prod_{n=1}^\infty p_n
\, \geq \, 1-\delta \, ,
\end{equation}
as desired.

We have thus far showed that $L=\infty$ almost surely. Take now an
ordered
countable dense deterministic subset $\{\theta_n\}_{n\geq1}$ of
$\{0\}\times[0,1]$
and let $\{\gamma_n\}_{n\geq1}$ be defined inductively as follows.
$\gamma_1$ is the subpath of the path from $\theta_1$ down to time $0$;
for $n\geq2$, $\gamma_n$ is the subpath of the path from $\theta_n$
down
to time either $0$ or when the latter path meets any of the $\gamma_i$,
$1\leq i\leq n-1$, whichever time is greater.
Then $\{\gamma_n\}_{n\geq1}$ is a disjoint family,
and the sum of the length of the $\gamma_n$'s, which equals $L$, is
almost surely infinite. This implies that there almost surely are
infinitely many marks in the union of the traces of the $\gamma_n$'s.
Since each $\gamma_n$ is finite, each one has almost surely finitely
many marks. This and the previous statement imply that there almost
surely are infinitely many marked $\gamma_n$'s and hence infinitely
 many
{\it distinct\/} nucleation points from among the
$\Theta(\theta_n)$'s.
Thus there are infinitely many colors for the $\theta_n$'s
in the $q=\infty$ case.

By the last part of property (3) of Theorem~\ref{teo:Thetatheorem}
 it follows that between every two $\theta_n$'s of different
color, there must occur at least one point on the interval with
(at least) those two colors, and so there are infinitely many
points with at least two colors. To see that these have zero
Lebesgue measure in the interval and that there are no points with
three (or more) colors, note that all such points must be double
points (i.e., two backward paths going down from that point) or
(if they have more than two colors) triple points of the backward
BW. Double points of the backward BW at $x=0$ correspond to
(non-starting point) zeros of the paths of the forward BW, which
have Hausdorff dimension $1/2$ and zero Lebesgue measure, while
triple points of the backward BW correspond to places where two
paths of the forward BW coalesce, which has zero probability of
occuring at a deterministic value $x=0$. This completes the proof
of Property (7).

   \subsection{Proof of Theorem~\ref{teo:Thetatheorem}
(8) -- color region Lebesgue measure }

By the sequential construction of Section~\ref{sec:mbw},
every nucleation point $(x',t')$ is the $j$'th marked point for some
$j$ along
the backward Brownian web path starting from some point
$(\bar{x}_i,\bar{t}_i)$ in a deterministic dense countable
set ${\cal D}$ of $\R^2$. Since the MBW distribution does not
depend on the ordering of ${\cal D}$, we will consider the
$j$'th mark on the path
from the first point $(\bar{x}_1,\bar{t}_1)$.
Furthermore since our arguments do not depend on the value of
$(\bar{x}_1,\bar{t}_1)$, we will take it to be the origin $(0,0)$.

Let $(x_k,\tau_k)$ denote the k'th marked point along the
backward web path starting at $(0,0)$. Our object is to prove
that for $k \geq 1$, the unique-color region,
\[
C_{(x_k,\tau_k)}^u = \{ (x,t): \Theta(x,t) = (x_k,\tau_k)\} \, ,
\]
has Lebesgue measure ${\cal L} (C_{(x_k,\tau_k)}^u) > 0$ a.s.

We first consider $k=1$.
Let $B_\varepsilon(x,t)  = (x-\varepsilon,x+\varepsilon)
\times(t-\varepsilon,t) $
and denote $B_\varepsilon (0,0)$ by $B_\varepsilon$.
We also define $Y_\varepsilon = {\cal L} (C_{(x_1,\tau_1)}^u \cap
B_\varepsilon)$ and
$X_\varepsilon = {(2 \varepsilon^2)}^{-1} Y_\varepsilon$ so that $0
\leq X_\varepsilon \leq 1$.
Since $Y_\varepsilon$ is decreasing in $\varepsilon$, to prove $\P(Y_1
> 0) =
1$, we will use that for $\varepsilon <1$,
\[
\P(Y_1 > 0) \geq \P(Y_\varepsilon > 0) \geq \P(Y_\varepsilon \geq
\varepsilon^2) =
\P(X_\varepsilon \geq 1/2) \, ,
\]
and argue that $X_\varepsilon \to 1$ in probability as $\varepsilon
\to 0$.
These imply that $\P(Y_1 > 0) \geq 1 - \delta$ for every $\delta > 0$
and thus
that $\P(Y_1 > 0) = 1 $ as desired. That $X_{\varepsilon} \to 1$
will be a consequence of showing that $E(X_\varepsilon) \to 1$ as
$\varepsilon
\to 0$ which we proceed to do now.

Using Fubini's Theorem, we have that
\[
\E (X_\varepsilon) = \frac{1}{2 \varepsilon^2} \int_{B_\varepsilon}
\P(\Theta (x,t) = \Theta
(0,0) ) dx dt
\]
Since $\psi(x,t) = \P(\Theta(x,t) = \Theta(0,0))$ is equal to
$1$ at $(x,t) = (0,0)$, to see that $\E(X_\varepsilon) \to 1$, it
suffices to show that $\psi$ is continuous at $(0,0)$. This
can be seen easily from ~\eqref{eq:psi1}, or can be shown
directly by considering two marked Brownian paths starting at
$(0,0)$ and $(x,t)$.

To extend the argument to $k \geq 2$, we note that the same
reasoning shows that it suffices to show that the expression
\begin{equation}\label{eq:color}
\frac{1}{2 \varepsilon^2} \int_{\R^2} \P [ 1\{(x,t) \in B_\varepsilon
(x_{k-1},\tau_{k-1})\} 1 \{(\Theta(x,t) =
\hat{\Theta}(x_{k-1},\tau_{k-1}))\} ] dx dt
\end{equation}
tends to $1$ as $\varepsilon \to 0$, where
$\hat{\Theta}(x_{k-1},\tau_{k-1})) = (x_k,\tau_k)$ is the
first mark strictly after $(x_{k-1},\tau_{k-1})$ along the
backward path from $(0,0)$.

In the above integral, $(x,t)$ is deterministic and by the
sequential construction of MBW, $\Theta(x,t)$ is simply the
first mark along the backward path from $(x,t)$. If we denote
the marked Brownian web (backward) path starting from $(0,0)$
by $\hat{B}_{(0,0)}(s)$ for $s \leq 0$, then by the strong Markov property
for (a single) marked Brownian motion we have
that ${\hat B}_{k-1}(s) :=\hat{B}_{(0,0)}(s+
\tau_{k-1}) - \hat{B}_{(0,0)}(\tau_{k-1}), s \leq 0$ is a standard
(reversed)
Brownian motion and ${\hat B}_{k-1}(s), s \leq 0$ is
independent of ${\hat B}_{(0,0)}(s), s \geq
\tau_{k-1}$.  Now taking conditional expectation with respect to
the value $(y,u)$ of $(x_{k-1},\tau_{k-1})$,
expression~\eqref{eq:color} becomes

\begin{eqnarray*}
& &\frac{1}{2 \varepsilon^2} \int_{\R^2} \E_{(y,u)}(\P [ 1\{(x,t) \in
B_\varepsilon
(y,u)\} 1 \{(\Theta(x,t) =
\hat{\Theta}(y,u))\}  
) dx dt \\
 & = & \frac{1}{2 \varepsilon^2} \E_{(y,u)}\int_{\R^2} \P [ 1\{(x,t)
\in B_\varepsilon
(y,u)\} 1 \{(\Theta(x,t) =
\hat{\Theta}(y,u))\}] dx dt\\
 &=&  \frac{1}{2 \varepsilon^2} \int_{B_\varepsilon} \P(\Theta (x,t) =
\Theta(0,0)) dx dt
\end{eqnarray*}
where we have used the independence of ${\hat B}_{k-1}(s), s \leq 0$
and ${\hat B}_{(0,0)}(s), s \geq
\tau_{k-1}$ in the second line.
This reduces the argument to the $k=1$ case, thus
proving the theorem.


\subsection{Proof of Theorem~\ref{teo:q=2}}
\label{subsec:finiteq}

We start by giving the relationship of the noisy voter model (NVM)
on $\Z$ and the stochastic Ising model (alluded to just before the
 statement of Theorem~\ref{teo:q=2}).

The stochastic Ising model at inverse temperature $\beta$ is an
 interacting particle system with state space $\{-1,+1\}^\Z$ whose
flip rate at $x\in\Z$ for a state $\sigma$ is given by
\cite{L1,kn:NSch}
\begin{equation}
 \label{eq:rates}
 c(x,\sigma)=\frac1{1+\exp\{\beta\sigma(x)[\sigma(x-1)+\sigma(x+1)]\}}.
\end{equation}
 One readily checks that this is equivalent to the NVM with
$p=\frac2{1+e^{2\beta}}$.
The invariant measure is the Ising model Gibbs measure for the
formal Hamiltonian~\cite{L1}
 \begin{equation}
 \label{eq:ham}
 H(\sigma)=-\frac12\sum_{x\in\Z}\sigma(x)\sigma(x+1)
 \end{equation}
at inverse temperature $\beta$. This is a stationary (spatial)
Markov chain with state space $\{-1,+1\}$, transition matrix
\cite{kn:G}
\begin{equation}
 \label{eq:mc}
 \frac1{1+e^\beta} \left(
 \begin{array}{cc}
   e^\beta & 1 \\
   1 & e^\beta
 \end{array}
 \right),
\end{equation}

and with single site uniform distribution on $\{-1,+1\}$. For this
chain, runs of $+1$'s and $-1$'s have i.i.d.~lengths with a common
geometric distribution of mean $e^\beta$. Under the rescaling in
the statement of Theorem~\ref{teo:Phitheorem}, the block lengths
are i.i.d.~geometrics with mean
$\delta^{-1}\sqrt{2/\lambda-\delta^2}$ multiplied by $\delta$.
They thus converge to i.i.d.~exponentials with mean
$\sqrt{2/\lambda}$. The limiting color configuration can be then
described as a stationary (spatial) Markov jump process with state
space $\{-1,+1\}$ and uniform jump rate $\sqrt{\lambda/2}$.

We want now to identify the fixed time color configuration (as a
function of $x$) of the two-color CNVM with the above jump
process. For that, we first note, by Property (6) of
Theorem~\ref{teo:Thetatheorem}, that the configuration can be
described as a $\{-1,+1\}$-valued jump process. To characterize
this process, it is thus enough to describe its finite dimensional
distributions. But, by Theorem ~\ref{teo:Phitheorem}, these are
limits of the scaled two-color NVM finite dimensional
distributions, which in turn, by the above paragraph, are the
finite dimensional distributions of the Markov jump process
described there.

\appendix


\section{Hausdorff dimension of the graph of the sum of two Brownian
motions}
\label{app:haus}
\setcounter{equation}{0}

\bprop
\label{prop:haus}
Let $X,Y$ be two independent standard Brownian motions and let ${\cal
T}^+$ denote the set
of record times of $X$, i.e., ${\cal T}^+=\{t\geq0:\,X(t)=M(t)\}$,
where $M(t):=\sup_{0\leq s\leq t}X(s)$ is the maximum of X up to time
$t$.
Then, for $t_0>0$ and $a,b\in\r$ with $|a|+|b|>0$, the set
${\cal G}^+:=\{(aX(t)+bY(t),t):\,t\in{\cal T}^+\cap[0,t_0]\}$,
the projection of ${\cal T}^+\cap[0,t_0]$ onto the graph of $aX+bY$,
has Hausdorff dimension $1$ almost surely.
\eprop

\noindent{\bf Proof } An upper bound of $1$ for the
Hausdorff dimension follows readily from the
fact that ${\cal G}^+$ is the image of a set, ${\cal T}^+\cap[0,t_0]$,
of Hausdorff dimension $1/2$ a.s.~(since ${\cal T}^+\cap[0,t_0]$ has
the same
distribution as the set of zeros of $X$, $\{t\in[0,t_0]:\,X(t)=0\}$;
this follows from $(M(t)-X(t):\,0\leq t\leq t_0)$ having the same
distribution
as $(|X(t)|:\,0\leq t\leq t_0)$~\cite{kn:L}) --- see~\cite{kn:Ta} ---
under a
map which is a.s.~(uniformly) H\"older continuous of exponent $\alpha$
for every
$\alpha<1/2$, namely, the map $t\to(aX(t)+bY(t),t)$, where we use the
well known H\"older continuity properties of Brownian motion.

The desired lower bound is obtained by noting that the Hausdorff
dimension of
${\cal G}^+$ is bounded below by the Hausdorff dimension of the image
of ${\cal T}^+\cap[0,t_0]$ under $aX+bY$, or equivalently under $aM+bY$,
namely
$\{aM(t)+bY(t):\,t\in{\cal T}^+\cap[0,t_0]\}$. Notice that the latter
set equals
$\{as+bY(T(s)):\,s\in[0,M(t_0)]\}$,
where $T$ is the hitting time process associated to $X$, defined as
$T(x) := \inf\{t\geq0:\,X(t)=x\}$.
It suffices to show that the Hausdorff dimension of
$\{as+bY(T(s)):\,s\in[0,L]\}$
is a.s.~greater than or equal to $1$ for every deterministic $L>0$.
But that follows from known results as well. $Z(t):=at+bY(T(t))$ is a
self similar process of
exponent $1$ with stationary increments and satisfies also the
following condition
of Theorem 3.3 in~\cite{kn:XL}, from which the dimension bound follows.
The condition is that
there exists a constant $K$ such that $\P(|Z(1)|\leq x)\leq Kx$ for
every $x\geq0$.
This property is readily obtained from the distributions of
$Y$ and the hitting time variable $T(1)$.


\noindent L.~R.~G.~Fontes \\
Instituto de Matem\'atica e Estat\'{\i}stica, Universidade de S\~ao
Paulo\\ Rua do Mat\~ao 1010, CEP 05508--090, S\~ao Paulo SP,
Brazil

\vspace{.5cm}

\noindent M.~Isopi\\
Dipartimento di Matematica, Universit\`a di Roma ``la Sapienza''\\
P.le Aldo Moro 5, 00185 Roma, Italia

\vspace{.5cm}

\noindent C.~M.~Newman\\
Courant Institute of Mathematical Sciences, New York University\\
New York, NY 10012, USA

\vspace{.5cm}

\noindent K.~Ravishankar\\
Department of Mathematics, SUNY-College at New Paltz \\
New Paltz, New York 12561, USA

\end{document}